%
%
\overfullrule=0pt
\baselineskip 14pt
\settabs 8 \columns
\+ \hfill&\hfill \hfill & \hfill \hfill   & \hfill \hfill  & \hfill \hfill   &
 & \hfill \hfill\cr
\parindent=.3truein
\hfuzz=3.44182pt
\hsize 6truein
\font\title=cmbx10 scaled\magstep2
\font\normal=cmr10 
\font\small=cmr6
\font\ninerm=cmr9 
\font\ita=cmti10
\font\bol=cmbx10
\font\normal=cmr10 
\font\ita=cmti10
\font\bol=cmbx10
\font\title=cmbx10 scaled\magstep2
\font\ninerm=cmr9 
\font\small=cmr6

\def \con {\subseteq}

\def\sig{\sigma}

\def \ggg {\gamma}

\def \-> {\rightarrow}

\def\DD {\Delta}

\def\om {\omega}
\def\la {\lambda}
\def\xon {x_1,x_2,\ldots ,x_n}

\def \sas {\vskip .06truein}
\def\sa{{\vskip .125truein}}

\def\sap{{\vskip .25truein}}

\def\aaa {\alpha}
\def\bbb {\beta}

\def\ggg {\gamma}
\def\aa {\alpha}

\def\gg {\gamma}
\def\con {\subseteq}
\def \ses {\enskip = \enskip}
\def \s=s {\enskip = \enskip}
\def \sps {\enskip + \enskip}
\def \sms {\enskip -\enskip}

\def \scs {\ssp , \ssp}
\def \ess {\enskip}
\def \ssp {\hskip .25em}
\def \bigsp {\hskip .5truein}
\def \part {\vdash}
\def \scos {\ssp :\ssp}
\normal

\vsize=8truein
\sap
\def\today{\ifcase\month\or
January\or February\or March\or April\or may\or June\or
July\or August\or September\or October\or November\or
December\fi
\space\number\day, \number\year}
\def\picture #1 by #2 (#3){
  \vbox to #2{
    \hrule width #1 height 0pt depth 0pt
    \vfill
    \special{picture #3} 
    }
  }

\font\small=cmr6
\def\CF {{\cal F}}

\def\RFX {Rat[{\cal F};x_1,\ldots ,x_n]}
\def \xsig {x_{\sig_1},x_{\sig_2},\ldots ,x_{\sig_n}}
\def \PP {\Phi}
\def \PS {\Psi}
\def \SYM {Sym[\CF;\xon]}
\def \RSYM {Ratsym[\CF;\xon]}
\def \part {\vdash}
\def \DD {\Delta}

\headline={\small
A. Garsia \hfill GALOISIAN  GALOIS THEORY $\ess\ess$\hfill\today $\ess\ess\ess$
\folio } \footline={\hfil}
\voffset=.5truein 
\centerline{\bol GALOISIAN GALOIS THEORY}
\sap
\centerline{\bol Lecture Notes in Computational  Algebra}
\sap

{\narrower\smallskip\noindent
{\bol ABSTRACT.}$\ess\ess\ess$ 
{\ninerm 
These notes are an exposition of Galois Theory from the original Lagrangian and
Galoisian point of view. A particular effort was made here to better understand the
connection between Lagrange's purely combinatorial approach and Galois algebraic
extensions of the latter. Moreover, stimulated by the necessities of present day
computer explorations, the algorithmic approach has been given 
priority  here over every other aspect of presentation. In particular, you may not find here the clean
simplistic look characteristic of the classical exposition of E. Artin. In contrast these notes
should provide a good starting point in attempting  constructions in this most
difficult computational arena.
}
\smallskip}
\sa

\noindent{\bol 1. Symmetric Functions}
\sa

Unless otherwise specified all fields we shall work with here will be assumed to have zero characteristic.
If $\CF$ is such a field and $\xon$ are indeterminates the expression 
$$
\CF[\xon ]
$$
will denote the ring of polynomials in $\xon$ with coefficients in $\CF$. In contrast to customary
notation, the field of rational functions of $\xon$ with coefficients in $\CF$ will be denoted
by 
$$
\RFX
$$
We may also write
$$
Rat[{\cal F};f_1,f_2,\ldots ,f_m] 
$$
to represent all rational expressions in $f_1,f_2,\ldots ,f_n$  with coefficients in $\CF$,
whatever $f_1,f_2,\ldots , f_m$ may be in any particular situation. 
\sa

As customary, $S_n$ denotes the group of all permutations of $1,2,\ldots ,n$.
\sas
 
\noindent
If $\sig= (\sig_1,\sig_2,\ldots, \sig_n)\in S_n$ and 
$$
\PP(\xon )\in \RFX
$$ 
we  set $\ess x\sig \ses \xsig$ and 
$$
\sig \PP \ses \PP(x\sig ) \ses \PP (\xsig )
$$
We shall say that $\PP(\xon )$ is $k-valued$ if and only if the collection
$$
\{\ssp \sig  \ssp \PP\ssp :\ssp \sig \in S_n \ssp \}
$$
has cardinality $k$. Of course $1$-valued functions are usually called {\ita symmetric}.
It will be convenient to denote here by $\SYM$ and  $\RSYM$ the collections of symmetric elements of 
$\CF[\xon ]$ and $\RFX$ respectively.  
\sa

The combinatorial study of $k$-valued function is properly the domain of Lagrange theory
and will be carried out in the next section. In this section we will limit ourselves to establishing
the results on symmetric functions that are needed in our further developments.
\sas

We recall that the symmetric polynomial
$$
e_k(\xon )\ses \sum_{1\leq i_1<i_2<\cdots <i_k\leq n}\ssp x_{i_1} x_{i_2}\cdots x_{i_k}\ess 
$$
is usually referred to as the $k^{th}\ess elementary$ symmetric function.

A vector of integers 
$$
\la \ses (\la_1 \geq \la_2\geq \cdots \geq \la_k>0)
$$
is said to be a partition of $n$ and we write $\la\part n$ if and only if
$$
\la_1+\la_2+\cdots +\la_k = n
$$
We shall also say that $k$ is {\ita the number of parts} of $\la$. 
If $\la$ and $\mu$ are partitions of the same number with $k$ and $h$ parts respectively,
we shall say that  $\la$ {\ita dominates} $\mu$ and write $\la \ssp \geq \mu$ if and only if
$$
\la_1+\la_2+\cdots +\la_s \geq \mu_1+\mu_2+\cdots +\mu_s\bigsp (\ess for \ess s=1,2,\ldots , min(h,k)\ess )
$$
It easy to see that this partial order is linearly extended by the lexicographic order of integer
vectors. Finally, given a partition $\mu$ the partition $\mu'$ whose parts $\mu_s'$ are
given by
$$
\mu_s'\ses \#\{i\ssp :\ssp \mu_i\geq s\}
$$
is usually referred to as the {\ita conjugate } of $\mu$. A simple combinatorial argument shows that
we have $\la\geq \mu$ if and only if $\la'\leq \mu'$. 
\sa

Given a vector
$$
p=(p_1,p_2,\ldots , p_n)
$$ 
of non-negative integers, the weakly decreasing rearrangement
of the positive components of $p$ will referred to as the {\ita shape} of $p$
and denoted by $\la(p)$. Given a partition $\la$ the symmetric polynomial
$$
m_\la (x) \ses \sum_{\la(p)=\la}\ssp x^p \bigsp (\ess x^p=x_1^{p_1}x_1^{p_2}\cdots x_n^{p_n}\ess )
$$
is usually referred to as the {\ita monomial symmetric function} indexed by $\la$. 
The following fact is immediate.
\sa

\noindent {\bol Theorem 1.1}
\sas

{\ita Every symmetric polynomial $\PP(\xon)$ has a unique expansion of the form
$$
\PP(x) \ses \sum_{\la} \ssp c_\la \ssp m_\la (x)
\eqno 1.1
$$
where the $c_\la$ are integers if $\PP$ has integer coefficients, and $c_\la$ is in $\CF$ if
$\PP\in \CF[\xon ]$. }
\sas
\noindent{\bol Proof}

The symmetry of $\PP$ implies that the coefficients $c_p$ and $c_q$ of any two monomials
$x^p$ and x$^q$ appearing in $\PP$ must be the same if $p$ and $q$ have the same shape.
Thus 1.1 is obtained by collecting terms of $\PP$ according to shape. 
\sa

Let now $\mu=(\mu_1,\mu_2,\ldots ,\mu_h)$ be a partition and set
$$
e_\mu (x)\ses e_{\mu_1}(x)e_{\mu_2}(x)\cdots e_{\mu_h}(x)
\eqno 1.2
$$
Since this polynomial is clearly symmetric in $\xon$ it must have an expansion of the form
$$
e_\mu(x) \ses \sum_\la\ssp m_\la(x) \ssp c_{\la\mu}\ess .
\eqno 1.3
$$
It develops that the coefficients $c_{\la\mu}$ have a suggestive combinatorial interpretation.
\sa

\noindent{\bol Proposition 1.1}
\sa

{\ita If $\mu$ has $h$ parts and $\la$ has $k$ parts then $c_{\la\mu}$ gives the number
of $h\times k$ matrices with $0,1$-entries and row and column sums given by $\mu$ and $\la$ respectively.}
\sas
\noindent{\bol Proof}

There is a one-to-one correspondence between these matrices and the monomials obtained by expanding
the product in 1.3. In fact, if $m_1(x),m_2(x),\ldots ,m_h(x)$ are monomials coming out of
$e_{\mu_1}(x),e_{\mu_2}(x),\ldots , e_{\mu_h}(x)$ respectively, then $m_i(x)$ corresponds to a subset
of $1,2,\ldots ,n$ of cardinality $\mu_i$, thus it may be represented by a $0,1$-vector with
$n$ components in an obvious manner. Putting together these vectors as the rows of an $h\times n$
matrix $M$ we see that these monomials multiply to $x^p=x_1^{p_1}x_2^{p_2}\cdots x_n^{p_n}$
if and only if the columns of $M$ add up to $p_1,p_2,\ldots ,p_n$ respectively. This is our desired 
correspondence. Since, since monomials of the same shape have the same coefficient, the assertion
follows by taking $x^p=x_1^{\la_1}x_2^{\la_2}\cdots x_k^{\la_k}$.
\sa

\noindent{\bol Remark 1.1}
\sas

Let $M$ be one of the $0,1$-matrices with row and column sums given by $\mu$ and $\la$. 
Note that the number of $1\ssp 's$ in the first $s$ columns of $M$ is equal to
$$
\la_1+\la_2+\cdots +\la_s
$$ 
On the other hand by moving all the $1\ssp 's$ along their rows until they are bumper to bumper to
the left (and all the $0's$ similarly to the right)
we obtain a matrix $M'$ whose column sums are the parts of the partition $\mu'$
conjugate to $\mu$. This not only gives the inequality
$$ 
\mu'_1+\mu'_2+\cdots +\mu'_s \ssp \geq \la_1+\la_2+\cdots +\la_s\bigsp (\ess for \ess s=1,2,\ldots ,min(h,k)
\ess )
$$
but also assures that when $\la=\mu$ there can only be one matrix with the desired row and column sums.
\sa

We can thus obtain
\sa

\noindent{\bol Theorem 1.2}
\sas

{\ita Every homogeneous symmetric polynomial $\PP(\xon)$ of degree $m$ has a unique expansion
of the form
$$
\PP(x) \ses \sum_{\mu\part m} \ssp d_\mu \ssp e_\mu (x)
$$
where the coefficients $d_\mu$ are integers if $\PP$ has integer coefficients and they
are elements of $\CF$ if $\PP\in \SYM$. In particular, $\SYM$ is the polynomial 
ring generated by the elementary symmetric functions $e_1(x),e_2(x),\ldots ,e_n(x)$.}
\sas

\noindent{\bol Proof}

Let $D=\|d_{\la\mu}\|_{\la,\mu\part m}$  denote the matrix obtained
when the partitions $\la$ are in lexicographic order and the partitions $\mu$ are in
reverse lexicographic order. This done we see  that an immediate consequence of the observations 
made in the Remark above is that $D$ must be unitriangular and therefore invertible over the integers.
This shows that the collection  $\{e_\mu(x)\}_{\mu\part m}$ must also give a basis and that
the elements of the basis  $\{m_\la(x)\}_{\la\part m}$ have integral linear expansions in terms
of the $e_\mu(x)\ssp 's$. This establishes our assertions.
\sa

\noindent {\bol Corollary 1.2}
\sas

$$
\RSYM \ses Rat[\CF;e_1,e_2,\ldots ,e_n]
\eqno 1.4
$$
\noindent{\bol Proof}

Note that every element $\PP\in \RFX$ can be expressed in the form
$$
\PP(x)={P(x)\over Q(x)}
$$
with $P,Q\in \CF[\xon]$. Multiplying numerator and denominator by the polynomial
$$
\prod_{\sig\in A_n \atop \sig\neq id}\ssp \sig Q
$$
we can write
$$
\PP(x)={P^*(x)\over Q^*(x)}
$$
where
$$
P^*(x)= P(x)\prod_{\sig\in A_n \atop \sig\neq id}\ssp \sig Q 
\ess\ess\ess\ess , \ess\ess\ess\ess 
Q^*(x) = \prod_{\sig\in A_n }\ssp \sig Q \ess . 
$$
Since by its very construction $Q^*(x)$ is already symmetric, we see that $\PP$ will be symmetric
if and only if $P^*$ is. Thus our assertion follows immediately from Theorem 1.2.
\vfill\supereject

\noindent{\bol 2.  The Euclidean algorithm the Resultant and the Discriminant.}
\sa

If $A,B$ are polynomials in $\CF[t]$, with $degree\ssp B<degree \ssp A$, 
then the {\ita quotient } and the {\ita remainder} anof the {\ita division} of $A$ by $B$ are
respectively the unique
polynomials $Q$ and  $R$ satisfying the requirements
$$
\eqalign{
&(1)\ess\ess A\ses B\ssp Q\sps R\cr
&(2)\ess\ess degree(R)\ssp<\ssp degree(B)\cr}\ess .
$$
To construct these two polynomials we can proceed as follows. Set
$$
\DD(A/B)\ses{leading\ess coeff \ess A\over leading\ess coeff  \ess B}\ess\ess  t^{degree(A)-degree(B)}
$$
Clearly we  have
$$
degree\bigl(A-\DD(A/B)\ssp B\bigr)\ess <\ess degree(A).
\eqno 2.1
$$
Now set
$$
R^{(o)}= A\ess\ess {\rm and} \ess\ess R^{(i)}= R^{(i-1)}- \DD(R^{(i-1)}/B)\ssp B 
\bigsp (\ess i=1,2,\ldots\ess )
\eqno 2.2
$$
Since, by 2.1, degrees are decreasing at least by one  at each step of the recursion, after
$k\leq degree(A)-degree(B)$ steps we shall have $degree(R^{(k)})<degree(B)$. At this point we stop
the recursion. By adding the identities in 2.1 we easily derive that
$$
A\ses \Bigl(\sum_{i=0}^{k-1}\ssp \DD(R^{(i)}/B)\ssp\Bigr)\ess B \sps R^{(k)}
\eqno 2.3
$$
Thus we may take
$$
Q= \sum_{i=0}^{k-1}\ssp \DD(R^{(i)}/B)
\ess\ess\ess , \ess\ess\ess R=   R^{(k)}\ess .
\eqno 2.4
$$
The important fact to note is that from 2.3 we deduce that if $A,B\in \CF[t]$ then
$R(t)$ is in $\CF[t]$ as well. We shall refer to the process above as the {\ita division algorithm}.
\sa

The Euclidean algorithm is the process which yields the greatest common divisor $D(t)$ 
of two polynomials $A$ and $B$. It is shown by  Berlekamp [1] that $D$ may be computed by the following
process. Set 
$$
r_{-2}=A
\ess\ess ,\ess\ess 
r_{-1}=B 
\ess\ess ,\ess\ess 
p_{-2}=0
\ess\ess ,\ess\ess
p_{-1}=1 
\ess\ess ,\ess\ess 
q_{-2}=1
\ess\ess ,\ess\ess 
q_{-1}=0
$$
Then compute $a_k,r_k,p_k,q_k$ according to the recursions
$$
\eqalign{
&(1)\ess\ess r_{k-2}\ses a_k\ssp r_{k-1}\sps r_k\ess\ess (division)\cr
&(2)\ess\ess p_k\ses a_k\ssp p_{k-1}\sps p_{k-2}\cr
&(3)\ess\ess q_k\ses a_k\ssp q_{k-1}\sps q_{k-2}\cr}\ess .
\bigsp (\ess k=0,1,2,\ldots\ess )
$$
since $degree(r_k)$ dcreases at least by one at each step after $n<degree(B)$ steps
we shall have $r_n=0$. It is shown in [] that these recursion force the following basic identities
$$
\eqalign{
&(1)\ess\ess q_n \ssp p_{n-1}\sms p_n\ssp q_{n-1}\ses =(-1)^n\cr
&(2)\ess\ess A\ses r_{n-1}\ssp p_n\cr
&(3)\ess\ess B\ses r_{n-1}\ssp q_n\cr}\ess .
$$
Since equation (1) here yields that $p_n$ and $q_n$ are relatively prime, we see that (2) and (3)
yield that the greatest common divisor $D$ of $A$ and $B$ is necessarily given by $r_{n-1}$.

As pointed out in [] the advantage of this process over the one that is usually described in 
most textbooks is that it provides the final answer without excessive storage of partial results.
In fact, only $7$ results need to be stored at any particular time, a number that is independent
of the choice of $A$ and $B$.
\sa

If the roots of a polynomial
$$
P(t)\ses a_o+a_1t +a_2t^2+\cdots +a_mt^m
$$
are $\xon$ then when $a_n=1$ we may write it in the form
$$
P(t)\ses (t-x_1)(t-x_2)\cdots (t-x_m)
\eqno 2.5
$$
Thus, we see that we must have
$$
a_{n-k}\ses (-1)^{k}\ssp e_k(\xon).
\eqno 2.6
$$
Clearly, if $P(t)$ and 
$$
Q(t)\ses b_o+a_1t +b_2t^2+\cdots +b_mt^m\ses (t-y_1)(t-y_2)\cdots (t-y_n)\bigsp (\ess b_n=1\ess )
$$
have a root in common, the expression
$$
\prod_{i=1}^m\prod_{j=1}^n(x_i-y_j)
$$
will necessarily vanish. It develops that a multiple of this expression may be written as a polynomial 
in the coefficients of $P$ and $Q$. This polynomial is usually referred to as the {\ita Resultant of $P$ 
and $Q$} and will be denoted here by $R[P,Q]$.
The case $m=2$, $n=3$ is sufficient to get accross the idea and avoids excessive notation.
\sa

\noindent{\bol Theorem 2.1}

$$
R[P,Q]= det \ess 
\pmatrix{
a_o & a_1 & a_2 & 0 & 0\cr
0 &a_o & a_1 & a_2 & 0 \cr
0 & 0 & a_o & a_1 & a_2 \cr
b_o & b_1 & b_2 & b_3 & 0\cr
0 &b_o & b_1 & b_2 & b_3 \cr
}
\ses (-1)^{2\times 3}\ssp a_2^3\ssp b_3^2\ess \prod_{i=1}^2  \prod_{j=1}^3(x_i-y_j)
\eqno 2.7
$$
{\bol Proof}

Note that since the division by $a_2$ and $b_3$ does not change the roots of $P$ and $Q$, we can
divide both sides of 2.7 by  $a_2^3\ssp b_3^2$ and reduce ourselves to the case $a_2= b_3=1$.
This given, note that we have the following matrix multiplication identity:
$$
\pmatrix{
a_o & a_1 & a_2 & 0 & 0\cr
0 &a_o & a_1 & a_2 & 0 \cr
0 & 0 & a_o & a_1 & a_2 \cr 
b_o & b_1 & b_2 & b_3 & 0\cr
0 &b_o & b_1 & b_2 & b_3 \cr
}
 \times 
\pmatrix{
1 & 1 & 1 & 0 & 0\cr
y_1 & y_2 & y_3 & 0 & 0\cr
y_1^2 & y_2^2 & y_3^2 & 0 & 0\cr
y_1^3 & y_2^3 & y_3^3 &1 & 0\cr
y_1^4 & y_2^4 & y_3^4 & 0 & 1\cr}
\ess =\ess
\pmatrix{
P(y_1) & P(y_2) & P(y_3) & 0 & 0\cr
y_1P(y_1) & y_2P(y_2) &y_3P(y_3) & 0 & 0\cr
y_1^2P(y_1) &y_2^2 P(y_2) & y_3^3P(y_3) & 0 & 0\cr
0 & 0 & 1 & b_3 & 0\cr
0 & 0 & 0 & 0 & b_3\cr}
$$
and 2.7 follows immediately by equating determinants of boths sides and cancelling
the common factor.  
\sa

\noindent{\bol Remark 2.1}
\sas

Note that the vanishing of the determinant in 2.7 assures that we can find a non trivial solution
to the corresponding homogeneous system. Now a simple computation shows that we have
$$
\pmatrix{ \aaa_o & \aaa_1& \aaa_2 &-\bbb_o&-\bbb_1}\ess \times\ess \pmatrix{
a_o & a_1 & a_2 & 0 & 0\cr
0 &a_o & a_1 & a_2 & 0 \cr
0 & 0 & a_o & a_1 & a_2 \cr 
b_o & b_1 & b_2 & b_3 & 0\cr
0 & b_o & b_1 & b_2 & b_3 \cr
}\ess  = 0
$$
if and only if
$$
( \aaa_o + \aaa_1t +\aaa_2t^2) \ssp P(t) \sms  (\bbb_o+\bbb_1t) \ssp Q(t) \ses 0
$$
and this is equivalent to the statement that $P$ and $Q$ have a non trivial common factor.
\sa

In the same vein we see that the expression
$$
\DD(x)\ses \prod_{1\leq i<j\leq n}(x_i-x_j)
$$
vanishes if and only if $P(x)$ as given by 2.5 has multiple roots. Since, $\DD(x)^2$
is clearly symmetric in $\xon$, Theorem 1.2 and 2.6 guarantee that the latter polynomial
should be expressible as a polynomial in the coefficients $a_o,a_1,\ldots ,a_n$ of
$P(t)$. In fact, we need only replace $Q$ by the derivative $P'(t)={d\over dt}P(t)$ 
in $R[P,Q]$ to obtain
\sa

\noindent{\bol Theorem 2.2}

$$
R[P,P'] \ses a_m^{m(m-1)}\ssp \prod_{i=1}^n\prod_{j=1}^n\ ^{(i)}\ssp (x_i-x_j)= (-1)^{n(n-1)/2}\ssp
a_m^{m(m-1)}\ssp \DD(x)^2
\eqno 2.8
$$
{\ita where the superscript ``$\ ^{(i)}$'' in the product is to indicate that the factor corresponding
to $j=i$ is to be omitted.} 

\noindent
{\bol Proof}

It is easily seen that the general form of 2.7 may also be written as
$$
R[P,Q]\ses \ssp (-1)^{nm}a_m^n \prod_{i=1}^n\ssp Q(x_1)Q(x_2)\cdots Q(x_n)
\ess 
\eqno 2.9
$$
Now we easily see that 
$$
P'(x_i) \ses a_n\ssp \prod_{j=1}^n\ ^{(i)}\ssp (x_i-x_j)
\eqno 2.10
$$
and 2.8 follows by setting $Q=P'$ in 2.9.
\sa

The polynomial $R[P,P'] $ is usually referred to as the {\ita discriminant} of the equation
$$
P(t)=0\ess .
$$
\sa

We shall terminante with a simple fact that will play a crucial role in the sequel. 
Recall that a polynomial $B(t)$ with coefficients in a field $\CF$ is said to be irreducible over $\CF$
if it does not admit a factorization
$$
B(t)=P(t)Q(t)
$$
into two polynomials $P,Q\in \CF[t]$ of strictly lesser degree than $B$. The Euclidean Algorithm immediately
yields that if $B(t)$ is irreducible and $A(t)$ is a polynomial in $\CF[t]$ which has a root
in common with $A(t)$ then $B(t)$ must be a factor of $A(t)$. The reason for this is that if
$A(t)$ and $B(t)$ share a root (a fact which may be verified by computing $R[A,B]$) then
the greatest common divisor $D(t)$ of $A$ and $B$  (as yielded by the process descrbed above), is
in $\CF[t]$ as well and we would be led to a contradiction unless $D$ is equal to a constant
multiple of $B$. The important conclusion we draw from this is that if a polynomial $A\in \CF[t]$
shares a root with an irreducible polynomial $B(t)\in \CF[t]$ then it must vanish for all the other
roots of $B(t)$. This fact has the following immediate extension.
\sa

\noindent{\bol Proposition 2.1}

{\ita Let $B(t)\in \CF[t]$ be irreducible in $\CF$ and let $\PP(t)\in Rat[\CF;t]$ vanish 
at one of the roots of $B(t)$ then $\PP(t)$ vanishes at all the other roots of $B(t)$.}
\sas
\noindent{\bol Proof}

By hypothesis $\PP(t)=P(t)/Q(t)$ with $P,Q\in \CF$ using the Euclidean Algorithm we can cancel out
(if necessary) the greatest common divisor of $P$ and $Q$ assure that $P$ and $Q$ have no common root.
But then $\PP(t)$ vanishes if and only if $P(t)$ does, and we are thus reduced to the case discussed above.
\vfill\supereject   
\def \om {\omega}
\def\CF {{\cal F}}

\def\RFX {Rat[{\cal F};x_1,\ldots ,x_n]}
\def \xsig {x_{\sig_1},x_{\sig_2},\ldots ,x_{\sig_n}}
\def \PP {\Phi}
\def \PS {\Psi}
\def \SYM {Sym[\CF;\xon]}
\def \RSYM {Ratsym[\CF;\xon]}
\def \part {\vdash}
\def \DD {\Delta}
\def \omi {\hskip -.01truein\ ^{(i)}}
\def \scos {\ssp :\ssp}
\def \bup {\hskip -.01truein}
\def \om {\omega}
\def \eon {e_1,e_2,\ldots ,e_n}
\def \SQR {{\root 2 \of R}}
\def \sh {{\scriptstyle{1\over 2}}}
\def \shf {{\scriptstyle{1\over 4}}}
\def \AO {{\root 3 \of {-\sh q+ \SQR}}}
\def \BO {{\root 3 \of {-\sh q- \SQR}}}

\headline={\small
A. Garsia \hfill 264C Lecture notes$\ess\ess$\hfill GALOIS THEORY $\ess\ess$\hfill\today $\ess\ess\ess$
\folio } \footline={\hfil}
\voffset=.5truein 

\noindent{\bol 3. The cubic and the quartic}
\sa

Formulas giving the general solution of the cubic equation
$$
E_3(t)\ses (t-x_1)(t-x_2)(t-x_3)\ses t^3-e_1t^2+e_2t-e_3\ses 0
\eqno 3.1
$$
where first discovered by Ferreo (sometimes before 1505) rediscoverd by Tartaglia 
and published by Cardano in 1545. Setting $\om=e^{{2\pi\ssp i}/3}$ the three roots of 3.1
may be written as follows

$$
\eqalign{
x_1 &\ses {e_1\over3}\sps \ess \ess\AO \ess\sps\ess  \BO\ess ,\cr
x_2 &\ses {e_1\over3}\sps \om\ssp \AO \sps\om^2\ssp \BO\ess ,\cr
x_3 &\ses {e_1\over3}\sps\om^2\ssp \AO \sps \om\ssp\BO\ess ,\cr}
\eqno 3.2
$$
where
$$
p\ses e_2-e_1^2/3
\ess\ess ,\ess\ess
q=-e_3+{{\scriptstyle{1\over 3}}}e_1e_2 -{\scriptstyle{2\over 27}}\ssp e_1^3
\ess\ess and \ess\ess R\ses (q/2)^2+(p/3)^3\ess .
\eqno 3.3
$$

These formulas are usually derived by the following process, apparently due to Hudde (1650).
We start by making the substitution $t=y+{\scriptstyle{1\over 3}}e_1$ in 3.1 and transform it to
$$
y^3+py+q=0
$$
This given the further substitution  
$$
y=z-{p\over3z}
\eqno 3.4
$$
brings us to  the equation
$$
z^3-{p^3\over 27z^3}+q\ses 0
$$
or better yet
$$
z^6+qz^3-{p^3/27}\ses 0\ess .
\eqno 3.5
$$
Since this is a quadratic equation for $z^3$ we immediately derive the two solutions
$$
z_1^3\ses -\sh\ssp q \sps \SQR
\ess\ess or \ess \ess 
z_2^3\ses -\sh\ssp q \sms \SQR
$$
since 
$$
(-\sh\ssp q \sps \SQR)(-\sh\ssp q \sms \SQR)\ses -(p/3)^3
$$
we may extract cube roots so that 
$$
z_2\ses\sms {p\over 3z_1}
\eqno 3.6
$$
This is given the six  roots of  3.5 are
$$
\matrix{
z_1&z_2\cr
\om z_1&\om^2 z_2\cr
\om^2 z_1& \om z_2\cr}
$$
Where each of these pairs multiplies to $-p/3$. 

We can now use 3.4 and derive that the three roots of 3.3 may be written in the form
$$
\eqalign{
y_1 &\ses \ess  z_1 \sps \ess  z_2\cr
y_2 &\ses \om z_1 \sps \om^2z_2\cr
y_3 &\ses \om^2 z_1 \ssp +\ssp \om z_2\cr}
\eqno 3.7
$$
From which the formulas in 3.2 can be immediately obtained.
\sa

The quartic equation was treated in a similar manner. That is ``ad hoc'' manipulations
were used to transform it to equations which could be solved by extraction
of roots. To give a brief idea of the process in this case we start with
$$
E_4(t)\ses (t-x_1)(t-x_2)(t-x_3)(t-x_4)\ses t^4-e_1t^3+e_2t^2-e_3t+e_4\ses 0\ess .
\eqno 3.8
$$
Completing the square suggested by the first two terms we can rewrite this equation in the form
$$
(t^2-\sh e_1\ssp t)^2\ses (\shf e_1^2-e_2)\ssp t^2 \sps e_3\ssp t  -e_4\ess .
$$
We then add $(t^2-\sh\ssp e_1t)y + {{\scriptstyle {1\over 4}}} y^2$ to both sides and get
$$
(t^2-\sh\ssp e_1t+\sh y )^2\ses 
({{\scriptstyle {1\over 4}}}e_1^2-e_2+y)\ssp t^2+(e_3-\sh e_1 y)\ssp t +
{{\scriptstyle {1\over 4}}}y^2 -e_4\ess .
\eqno 3.9
$$
Next,  $y$ is determined so that also the term on the right becomes a perfect square.
This requires the coefficients
$$
A= \shf e_1^2-e_2 +y  \ess\ess\ess ,\ess\ess\ess 
B= e_3-\sh\ssp e_1\ssp y \ess\ess\ess ,\ess\ess\ess 
C= \shf \ssp y^2-  e_4
$$
of the quadratic on the right hand side of 3.9 satisfy the equation 
$$
B^2-4AC\ses 0\ess .
$$
This leads to a cubic equation for $y$. To see what are the roots of this equation, 
we should try to factor $B^2-4AC$. Nowdays this is easily done, using any of 
the available computer algebra packages. In this manner we discover the pleasing
fact that $B^2-4AC$ factors beautifully in  terms of the roots of $E_4(t)$. Namely, we have
$$
B^2-4AC\ses  ( y - x_2 x_4 - x_1 x_3)\ess ( y- x_1 x_2 - x_3 x_4)\ess ( y- x_1 x_4 - x_2 x_3)\ess .
$$
Thus the three roots of this cubic when expressed in terms
of the roots of 3.8 are none other than
$$
y_1 \ses x_1x_2+x_3x_4
\ess\ess ,\ess\ess 
y_2 \ses x_1x_3+x_2x_4
\ess\ess ,\ess\ess 
y_3 \ses x_1x_4+x_2x_3
$$
We should also note that setting $y=y_1$ we have
$$
\eqalign{
&A\ses \shf e_1^2-e_2 +y_1 \ses \shf\ssp (x_1+x_2-x_3-x_4)^2\cr
&B\ses e_3-\sh\ssp e_1\ssp y_1\ses -\sh\ssp (x_1x_2-x_3x_4)(x_1+x_2-x_3-x_4)\cr 
&C\ses \shf \ssp y_1^2-e_4\ses \shf\ssp (x_1x_2-x_3x_4)^2\cr}
\eqno 3.10
$$
which yield that the right hand side of 3.9 when $y=y_1$ reduces to
$$
\shf\ssp \bigl(\ssp -(x_1+x_2-x_3-x_4)\ssp t +x_1x_2-x_3x_4\ssp \bigr)^2\ess .
$$
This allows us to rewrite 3.9 in the form
$$
L^2 - R^2\ses (L\sms R)\ess (L\sps R)\ses 0\ess ,
\eqno 3.11
$$
with
$$
L\ses t^2-\sh\ssp e_1t+\sh y_1 
\ess\ess\ess {\rm and}\ess\ess\ess 
R \ses \sh\ssp ( -(x_1+x_2-x_3-x_4)\ssp t +x_1x_2-x_3x_4) 
$$
Now the factorization of $A$ in 3.10 suggests setting
$$
z\ses x_1+x_2-x_3-x_4\ses {\root 2 \of {e_1^2-4e_2 +4y_1}}\ess ,
\eqno 3.12
$$
and then the factorization of $B$ in 3.10 gives that
$$
-\sh(x_1x_2-x_3x_4)\ses ( e_3-\sh e_1 y_1)/z\ess .
$$
This given, we may write
$$
-R\ses\sh z \ssp t \sps (e_3-\sh\ssp e_1y_1)/z
$$ 
and the equation in 3.11 may yet be rewritten as
$$
\bigl(t^2-\sh ( e_1+z)\ssp t+\sh y_1 -( e_3-\sh e_1 y_1)/z\bigr)\ssp
\bigl(t^2-\sh ( e_1-z)\ssp t+ \sh y_1 +( e_3-\sh e_1 y_1)/z\bigr)= 0\ess .
$$
Now this is none other than factoring $E_4(t)$
in the form
$$
E_4(t)\ses \bigl(t^2-(x_1+x_2)\ssp t -x_1x_2\bigr)\ess \bigl(t^2-(x_3+x_4)\ssp t -x_3x_4\bigr)
\ess .
$$
In fact, it may be easily verified that
$$
\eqalign{
&x_1+ x_2 \ses \sh e_1 +\sh z
\ess\ess ,\ess\ess 
x_1x_2\ses \sh y_1 -( e_3-\sh e_1 y_1)/z\ess ,\cr
&x_3+ x_4 \ses \sh e_1 -\sh z
\ess\ess ,\ess\ess 
x_3x_4\ses \sh y_1 +( e_3-\sh e_1 y_1)/z\ess .\cr}
$$
Thus we may obtain the desired expressions for the pairs $x_1,x_2$ and $x_3,x_4$  by solving 
the two quadratic equations
$$
\eqalign{
&t^2-\sh ( e_1+z)\ssp t+ \sh y_1 -( e_3-\sh e_1 y_1)/z\ses 0\ess ,\cr
&t^2-\sh ( e_1-z)\ssp t+ \sh y_1 +( e_3-\sh e_1 y_1)/z\ses 0 \ess . \cr }
$$
\sa

This given, in the $17^{th}$ and $18^{th}$ centuries it was natural to assume that
the solution of the general polynomial equation should be obtainable by similar
manipulations and successive root extractions. This was the motivating force []
in Lagrange's investigations in the 1770's that led him to his historic paper
{\ita R\'eflections sur la r\'esolution alg\'ebrique des \'equations }.

To be precise Lagrange was investigating the possibility of finding {\ita closed form}
expressions for the roots $\xon$ which (like those appearing in 3.2) only involved the
elementary symmetric functions $\eon$, roots of unity and radicals. We shall refer to this
as {\ita ``solving the general equation by radicals''}.

His point of departure was a close examination of the solutions of the cubic and the quartic.
Remarkably, he was able to sort out of those simingly ad hoc manipulations a unifying
general mechanism of solution. As we shall see Lagrange discovered that in both cases
the final formulas could be reached by a sequence of identical, {\bf purely combinatorial,} steps.

This done, he tried to apply this mechanism to the quintic only to discover
that the possibility of pushing it through to the production of general formulas
for the roots of the quintic appeared to lead to a contradiction! 

In fact, he was (and he knew he was) within reach  of proving the impossibility 
of solving the quintic equation by radicals. 

Around 1799 Ruffini tried to complete Lagranges proof and
although he was able to push the argument quite a bit further he nevertheless
was left with a hypothesis which he could not remove. 

The glory of proving the impossibility of solving the general equation by radicals
was bestowed to Abel (for his 1826 paper) (see []) even though he was only concerned with the quintic and, 
as in Ruffini's work, there were still a number of gaps in his arguments. We shall not
deal with Abel's work here since it it departs from the combinatorial approach proposed by Lagrange
and later completed by Galois. In fact, the missing step needed to complete Lagrange
argument and obtain the unsolvability of the quintic by radicals  can be supplied by one single idea
of Galois. 

\sa

To appreciate the beauty of Lagrange's discoveries we should view his results in the original 1771 form. 
Unfortunately, for clarity we must deviate a bit from Lagrange's terminology.
For instance, although Lagrange proved that the order of a subgroup of a group
is a divisor of the order of the group, he had to do so in an indirect manner,
since the notion of a {\ita group} in its present form really started with Galois. 
Although using modern terminology distorts somewhat the historical perspective,
we will try as much as possible  to keep unchanged the contents of Lagrange's discoveries. 
Our main goal in the next two sections is to present the basic theorems of what is now
referred to as {\ita Galois} theory in a sequence that makes the transition from 
Lagrange to Galois as natural and effortless as possible.
\sap

\noindent{\bol 4. Lagrange's ``Galois'' Theory}
\sa

Throughout  Lagrange's work the roots of an equation
$$
E_n(t)\ses a_o+a_1 t + a_2 t^2+\cdots +a_nt^n\ses =0
$$
are assumed to be independent variables  $\xon$ and $E_n(t)$ is written in the form
$$
E_n(t)\ses (t-x_1)(t-x_2)\cdots (t-x_n)
\eqno 4.1
$$
The basic idea that led Lagrange to an {\ita understanding } of the classical solutions of the
cubic and the quartic is a careful analysis of the effect that permutations of the roots
have on various rational functions of the roots. 
To make precise what we mean by this we need some notation.

We are given a field $\CF$ which remains unchanged throughout, and for a function $\PP$ of the roots which
may be in $\CF[\xon]$ or in $\RFX$ as needed,  we set 
$$
G_\PP \ses \{\ssp \sig\in S_n \scos \sig \PP=\PP\ssp\}\ess .
\eqno 4.2
$$
Although Lagrange did not realize (nor did he need) that $G_\PP$ is a group, we shall not
ignore this fact here and obtain Lagrange's results by standard present day techniques.
We recall that $G_\PP$ is usually referred to as the {\ita stabilizer} of $\PP$.
\sa

The first basic result of Lagrange can be stated as follows
\sas

\noindent{\bol Theorem 4.1}
\sa

{\ita For $\PP,\PS\in \RFX$ we have 
$$
G_\PS\con G_\PP
\eqno 4.3
$$
If and only if} 
$$
\PP\in Rat[\CF,e_1,e_2,\ldots  , e_n,\PS]
$$
{\bol Proof}

If 
$$
\theta (y_1,y_2,\ldots ,y_n,t)\ess \in\ess  Rat[\CF,y_1,y_2,\ldots  , y_n,t] 
$$ 
and 
$$
\PP(\xon) \ses \theta (e_1,e_2,\ldots ,e_n,\PS)\ess .
\eqno 4.4
$$
Then we clearly have 4.3 since every permutation $\sig$ leaves $e_1,e_2,\ldots ,e_n$ unchanged 
and if $\sig\in G_{\PS}$ then also $\PP$ does not change. So the condition in 4.3 is trivially necessary.

To show the converse, we resort to the left coset decomposition
$$
S_n\ses \tau_1\ssp G_\PS \sps \tau_2\ssp G_\PS \ssp +\cdots +\ssp \tau_k\ssp G_\PS \ess .
\bigsp (\ess \tau_1=identity\ess )
\eqno 4.5
$$
(which by the way, was a Lagrange invention) and set
$$
Q(t)\ses \sum_{i=1}^k\ssp \tau_i \PP \ssp \prod_{j=1}^k\omi\ssp (t-\tau_j\PS)
\ess .
$$
Since any $\sig\in S_n$ permutes the left cosets of $G_\PS$ we may write 
$$
\sig\ssp \tau_i\ses \tau_{\pi_i}h_i\bigsp (\ssp h_i\in G_\PS\ssp )
$$
where the map $i\rightarrow \pi_i$ is a permutation of $(1,2,\ldots ,k)$. In particular, from
4.3 we deduce that $\sig\ssp \PP=\tau_{\pi_i}\PP$, and thus we must have
$$
\sig Q(t)\ses  \sum_{i=1}^k\ssp \tau_{\pi_i} \PP \ssp \prod_{j=1}^k
\bup \ ^{(\pi_i)} \ssp (t-\tau_j\PS)\ses Q(t)\ess .
$$
This implies that the coefficients of $Q(t)$ are in $\RSYM$, so by Theorem 1.2 they 
are in $Rat[\CF,\eon]$. The same can be said about the polynomial
$$
P(t)\ses \prod_{i=1}^k(t-\tau_i\PS)\ess .
\eqno 4.6
$$
Now setting $t=\PS$ in $Q(t)$ we get
$$
Q(\PS)\ses \PP\ssp \prod_{j=1}^k\bup \ ^{(1)}\ssp (\PS-\tau_j\PS)
\ses \PP\ssp P'(\PS)\ess .
\eqno 4.7
$$
Since, by construction, the values $\tau_i\PS$ (for $i=1,2,\ldots ,k$) are all distinct we shall have 
$ P'(\PS)\neq 0$ and we can divide it out in 4.7 to obtain
$$
\PP\ses { Q(\PS)\over P'(\PS)} \ses \theta(\PS)
$$
with 
$$
\theta (t)\ses  {Q(t)\over P'(t)}\ess \in \ess Rat[\CF;\eon,t]
$$
as desired.
\sa

\def \pac {\ess}

\noindent{\bol Remark 4.1}

If $G$ is a group and $H\con G$ is a subgroup then the left coset decomposition
$$
G\ses \tau_1\ssp H\sps \tau_2\ssp H\sps\cdots\sps \tau_m\ssp H \bigsp (\ess \tau_1=identity\ess )
\eqno  4.8
$$
yields that 
\sas

\centerline{\ita If $H\con G$ then the order of $H$ divides the order of $G$}  \hfill 4.9
\sas

In fact, the equation above gives that $|G|/|H|=m$. This result, which appeared for the first time 
in the work of Lagrange, was formulated  there as a statement concerning the {\ita number} of different values
taken by rational functions of the roots. More precisely, we can derive from 4.5  that $\PS$ is a $k$-valued
rational function of the roots if and only if $|S_n|/|G_\PS|=k$. Similarly, $\PP$ is $ h$-valued if and only
if  $|S_n|/|G_\PP|=h$. This gives that when $G_\PS\con G_\PP$ we have $k/h=|G_\PP|/|G_\PS |$.
We shall refer to $h$ and $k$ respectively as the {\ita multiplicities} of $\PS$ and $\PP$. 
So taking $H=G_\PS$ and $G=G_\PP$ in 4.8 we get that $k=hm$. We see then that in Lagrange's language 
the statement in 4.9 becomes
$$
If\pac  \PP\pac is\pac a\pac rational \pac function\pac of\pac \PS\pac then\pac
 the\pac multiplicity\pac of\pac\PP\pac divides 
\pac that\pac of\pac\PS\pac .
$$

Lagrange's proof was based precisely on the coset decomposition. Only he did not have to call it that way.
Indeed to get 4.8 for $H=G_\PS$ and $G=G_\PP$ all he had to do was {\ita bunch together} the elements
of $G_\PP$ that yielded the same value of $\PS$. 
\sa

It will be good here and after for $H$ a subgroup of $G$ to express the fact that $m=|G|/|H|$
by writing $H\con_mG$. This given, these observations can be sharpened into the following corollary 
of Theorem 4.1. 
\sa

\noindent{\bol Theorem 4.2}
\sa

{\ita If for $\PP,\PS \in \RFX$  we have 
$$
G_\PS\con_k G_\PP
\eqno 4.10
$$
then $\PS$ satisfies an equation of degree $k$ with coefficients in $Rat[\CF;\eon,\PP]$
which is irreducible in $Rat[\CF;\eon,\PP]$.}
\sas

\noindent{\bol Proof}

Let 
$$
G_\PP\ses 
\tau_1\ssp G_\PS\sps \tau_2\ssp G_\PS\sps\cdots\sps \tau_k\ssp G_\PS \bigsp (\ess \tau_1=identity\ess )
\eqno  4.11
$$
and set
$$
Q(t)\ses \prod_{i=1}^k \ssp (t-\tau_i\ssp \PS)\ses q_o(x)+q_1(x)\ssp t+\cdots +q_k(x)\ssp t^k \ess .
\eqno  4.12
$$
Since any $\sig\in G_\PP$ permutes the left cosets of $G_\PS$ in 4.11 we may write 
$$
\sig\ssp \tau_i\ses \tau_{\pi_i}h_i\bigsp (\ssp h_i\in G_\PS\ssp )
$$
where the map $i\rightarrow \pi_i$ is again a permutation of $(1,2,\ldots ,k)$.
This gives that for all $\sig\in G_\PP$ we have
$$
\sig Q(t)\ses \prod_{i=1}^n\ssp (t- \tau_{\pi_i}\PS)\ses Q(t)\ess .
$$
Consequently each of the coefficients $q_i(x)$ is left unchanged by the elements of $G_\PP$.
From Theorem 4.1 we then derive that each $q_i(x) \in Rat[\CF,\eon,\PP]$. On the other
hand from 4.12 we get that
$$
Q(\PS) \ses  q_o(x)+q_1(x)\PS+\cdots +q_k(x)\PS^k \ses 0\ess .
$$
Now suppose, if possible, that $Q(t)$ has a factorisation $Q(t)=Q_1(t)Q_2(t)$ where both polynomials
$Q_1(t)$ and $Q_2(t)$ have coefficients in $Rat[\CF;\eon,\PP]$ It will then follow that
both of them will be invariant under the action of $G_\PP$. In particular,
we must have
$$
\tau_i(Q_1(\PS))\ses Q_1(\tau_i\PS)\bigsp (\ess  for\ess  any\ess\ess i=1,2,\ldots ,k\ess )
$$
So if $\PS$ is a root of the equation $Q_1(t)=0$ then all the other roots of $Q(t)$ must 
satisfy it as well and $Q_2(t)$ must reduce to a constant in $Rat[\CF,\eon,\PP]$. The
analogous conclusion holds if $Q_2(\PS)=0$. Thus $Q(t)$ is irreducible as asserted.
This completes our proof.
\sa

\vbox{
\noindent{\bol Remark 4.2}

Here and after, if $G_\PS\con _k G_\PP$ and we have the left coset decomposition in 4.11,
then functions 
$$
\PS_1= \tau_1\PS \scs \PS_2=\tau_2\PS\scs \ldots , \PS_k= \tau_k\PS
$$
will be referred to as the {\ita conjugates} of $\PS$ in $G_\PP$. Note that if $\PP$
is in $Rat[\CF,\eon,\PS]$, Theorem 4.1 assures that $G_\PS\in G_\PP$. So in any case we must 
have 4.11 for some $k$. Now suppose that $\PS$ is a root of the equation $R(t)=0$ where $R(t)$ 
is a polynomial of  degree $h$ with coefficients in $Rat[\CF,\eon,\PP]$. Since 
this polynomial is then invariant under the action of $G_\PP$  all the conjugates of $\PS$
in $G_\PP$ must also be roots of $R(t)=0$. This implies that the polynomial $Q(t)$ in 4.12 must 
be a factor of $R(t)$.  However, if $R(t)$ is also irreducible in $Rat[\CF,\eon,\PP]$, 
then $R(t)$ and $Q(t)$ can only differ by a factor in $Rat[\CF,\eon,\PP]$ and we must also $h=k$. 
This should explain why we call $PS_1,\PS_2,\ldots ,\PS_k$ the ``conjugates'' of $\PS$.
} 
\sa

We have reached a point where to proceed further we need to make more precise
what we mean by {\ita solving the general equation by radicals}. To begin with we
shall assume that the given field $\CF$ (nowdays referred to as the {\ita ground field}) 
contains all the roots of unity of any order $\leq n$. This given, 
{\ita solving by radicals} the $n^{tic}$ in 4.1, {\ita in the Lagrange setting}
is to mean that we can find a sequence of rational functions 
$\PP_i\in Rat[\CF;\xon]$ ($i=0,1,2,\ldots $) such that
$$
\eqalign{ 
\PP_o &\ssp\in\ssp Rat[\CF;\eon]\cr
\PP_i &\ses {\root  p_i\of {\theta_i(\eon,\PP_{i-1})}}\cr}
\eqno 4.13
$$
where each $\theta_i$ is a rational function
$$
\theta_i\ses \theta_i(y_1,y_2,\ldots, y_n,t)\ess\in\ess  Rat[\CF,y_1,y_2,\ldots, y_n,t]\ess .
\eqno 4.14
$$
Finally, we shall require that the end function of this sequence say $\PP_d$ be one of the roots or
better yet (as we shall see) a function from which all the roots may be derived by rational
operations. 
\sas

Note first that since for any two integers $p$ and $q$ and for any $\PP$ we have
$$
{\root pq \of {\PP}}\ses {\root p  \of {\root q \of {\PP}}}
$$ 
there is no loss in requiring that the integers $p_i$ in 4.13 are all primes. 
Finally we can simplify the convoluted form of the recursion in 4.14 
by rewriting it in the form
$$
\eqalign{ 
&a)\ess\ess \PP_o\ess  \in Rat[\CF;\eon]\cr
&b)\ess\ess\PP_i^{p_i} \in Rat[\CF,\eon,\PP_{i-1}) \cr}
\eqno 4.15
$$
\sas

\noindent{\bol Remark 4.3}
\sas

Note further that as long as the coefficients $m_1,m_2,\ldots ,m_n$ are all 
distinct the function
$$
v(x)\ses v(\xon)\ses m_1 x_1+m_2 x_2+\cdots + m_n x_n
\eqno 4.16
$$
will necessarily be $n!$-valued. Since its stabilizer consists of just the identity permutation,
the hypotheses of Theorem 4.1 are satisfied for any rational function of the 
roots $\xon$. Thus for any $\PP\in \RFX$ we can construct a rational function
$\theta_\PP(t)\in Rat[\CF;\eon ,t]$ giving
$$
\PP(\xon) \ses \theta_\PP\bigl( v(\xon)\bigr)\ess .
\eqno 4.17
$$
We should note here for further reference that the proof of Theorem 4.1 yields that
$$
\theta_\PP \ses {Q_\PP(t)\over P'(t)}
\eqno 4.18
$$
where 
$$
Q_\PP(t)\ses \sum_{\sig\in S_n}\ssp \sig\PP \prod_{\tau\in S_n}\bup ^{(\sig)}\ssp (t-\tau v)
\eqno 4.19
$$
and
$$
P(t)\ses \prod_{\tau\in S_n} (t-\tau v)\ess .
\eqno 4.20
$$
It is important to notice that the denominator of $\theta_\PP$ in 4.18 is independent of
$\PP$ itself.
\sa
\def \om {\omega}
\def\CF {{\cal F}}

\def\RFX {Rat[{\cal F};x_1,\ldots ,x_n]}
\def \xsig {x_{\sig_1},x_{\sig_2},\ldots ,x_{\sig_n}}
\def \PP {\Phi}
\def \PS {\Psi}
\def \SYM {Sym[\CF;\xon]}
\def \RSYM {Ratsym[\CF;\xon]}
\def \part {\vdash}
\def \DD {\Delta}
\def \omi {\hskip -.01truein\ ^{(i)}}
\def \scos {\ssp :\ssp}
\def \bup {\hskip -.01truein}
\def \om {\omega}
\def \eon {e_1,e_2,\ldots ,e_n}
\def \SQR {{\root 2 \of R}}
\def \sh {{\scriptstyle{1\over 2}}}
\def \shf {{\scriptstyle{1\over 4}}}
\def \AO {{\root 3 \of {-\sh q+ \SQR}}}
\def \BO {{\root 3 \of {-\sh q- \SQR}}}
\def \nor {\triangleleft}
\def \son {\supseteq}
\def \TT  {\Theta}
\sa
\noindent{\bol 5. Lagrange's derivation of the roots of the cubic and the quartic.} 
\sa

Armed with this information, Lagrange could then come up with the following a-priori reconstruction
of the solutions of the cubic and the quartic.
\sa

\noindent
a) {\ita The cubic}
\sas

Let $\om=e^{2\pi i/3}$ and note that, since $1,\om , \om^2$ are distinct, the expression
$$
(x_1+\om x_2+\om^2 x_3)/3
$$
is necessarily a $6$-valued function of the roots of $E_3(t)$. Its values are
\sas
$$
\matrix{
z_1=(x_1+\om x_2+\om^2 x_3)/3\ess\ess & z_2=(x_1+\om x_3  + \om^2 x_2)/3 \ess\ess\cr
z_3=(x_2+\om x_3+\om^2 x_1)/3\ess\ess & z_4=(x_2+\om x_1  + \om^2 x_3)/3 \ess\ess\cr
z_5=(x_3+\om x_1+\om^2 x_2)/3\ess\ess & z_6=(x_3+\om x_2  + \om^2 x_1)/3 \ess\ess\cr}\ess .
\eqno 5.1
$$
\sas

\noindent
Now we see that
$$
z_3=\om^2 z_1\ess \scs\ess z_5=\om z_1\ess\ess\ess {\rm and }\ess\ess\ess z_4=\om z_2\ess\scs\ess z_6=\om^2
z_2  \ess .
$$
Thus 
$$
\eqalign{
(z-z_1)(z-z_3)(z-z_5)&\ses z^3-z_1^3\cr
(z-z_2)(z-z_4)(z-z_6)&\ses z^3-z_2^3\cr
}
$$
This implies that the $6$-degree equation
$$
(z-z_1)(z-z_2)(z-z_3)(z-z_4)(z-z_5)(z-z_6)\ses 0
\eqno 5.2
$$
must take the form
$$
(z^3-z_1^3)(z^3-z_2^3)\ses z^6-(z_1^3+z_2^3)\ssp z^3 \sps z_1^3z_2^3\ses 0\ess .
\eqno 5.3
$$
Lagrange, on the basis that 5.2 is symmetric in the roots of $E_3(t)$,
could now predict that  the expressions
$$
z_1^3+z_2^3\ess\ess {\rm and} \ess\ess z_1^3z_2^3
$$
must necessarily be polynomials in $e_1,e_2,e_3$. And indeed it can directly be checked from 
5.1 that
$$
z_1^3+z_2^3\ses -q \ess\ess {\rm and} \ess\ess z_1^3z_2^3\ses -(p/3)^3
\eqno 5.4
$$
with $p$ and $q$ given by 3.3. This immediately leads to the bicubic equation
$$
z^6+q z^3 -(p/3)^3\ses 0 \ess .
\eqno 5.5
$$
Now its solution leads to the extraction of the square root of
$$
\PP_o= q^2 + 4(p/3)^3 = 4\ssp R \ess .
$$
This expression is essentially the discriminant of $E_3(t)$. More precisely we have
$$
\PP_o(x) \ses -{\scriptstyle{1\over 27}}(x_1-x_2)^2(x_1-x_3)^2(x_2-x_3)^2 
$$
and its square root may be chosen to be
\def \sqt {{\root 2 \of 3}}
$$
\PP_1(x)\ses i({x_1-x_2\over \sqt})({x_1-x_3\over \sqt})({x_2-x_3\over \sqt})\ess . 
$$
This is a polynomial in the roots of $E_3(t)$ whose stabilizer $G_{\PP_1}$ is the group
of even permutations of $S_3=G_{\PP_o}$. Now this is in perfect agreement with Theorem 4.2.
\sas

Of course, in view of 5.4, we also have
$$
\PP_o(x) \ses (z_1^3 +z_2^3)^2-4(z_1^3z_2^3)\ses (z_1^3-z_2^3)^2
$$
and it can be easily verified that
$$
\PP_1(x)\ses z_1^3-z_2^3\ess .
$$
In other words $\PP_1(x)$ is the solution of
$$
\PP_1^2(x)\ses \PP_o\ess .
$$
Using 5.4 again we deduce that
$$
z_1^3\ses -\sh q \sps \sh \PP_1(x)\ess .
$$
This given, the final step is the construction of the $6$-valued function
$$
\PP_2(x)= z_1=  (x_1+\om x_2+\om^2 x_3)/3
$$
whose stabilizer $G_{\PP_2}$ is trivial and is the solution of
$$
\PP_2(x)^3\ses -\sh\ssp q\sps \sh \PP_1(x)\ess .
$$
In summary, this construction of a $3!$-valued function of the roots of $E_3(t)$ has led us 
to the following scheme:
$$
\matrix{
function& & expression & group\cr
\ess& \ess  & \ess & \ess\cr
\PP_o& \ses &-{\scriptstyle{1\over 27}}(x_1-x_2)^2(x_1-x_3)^2(x_2-x_3)^2 & S_3\cr  
\downarrow& \ess &\ess &3\downarrow\cr
\PP_1&\ses & i({x_1-x_2})({x_1-x_3})({x_2-x_3})/3\sqt& \{id,(1,2,3),(1,3,2)\}\cr 
\downarrow& \ess &\ess &2\downarrow\cr
\PP_2&\ses &(x_1+\om x_2+\om^2 x_3)/3& \{id\}\cr}
$$
Moreover we have that
$$ 
\PP_2^3\ses -\sh\ssp q +\sh \ssp \PP_1
\ess\ess\ess ,\ess\ess\ess 
\PP_1^2\ses \PP_o
$$
This shows that the solution of the cubic can be obtained the succession of steps 
$$
\matrix{
\PP_2^3\in  Rat[\CF,e_1,e_2,e_3,\PP_1]& \longleftarrow &\PP_1^2\in  Rat[\CF,e_1,e_2,e_3,\PP_o]
& \longleftarrow & \PP_o\in  Rat[\CF,e_1,e_2,e_3]\cr
G_{\PP_2}\con G_{\PP_1}&  &  G_{\PP_1}\con G_{\PP_o} &  & G_{\PP_o}=S_3 \cr       }
$$
\sap

\noindent
b)  {\ita The quartic}
\sas

We can proceed in the same manner as for the cubic and construct a sequence of rational functions 
of the roots according to the scheme expressed in 4.15, terminating again with a $4!$-valued function  
$$
v(x)\ses m_1x_1+ m_2x_2+ m_3x_3+ m_4x_4\ess .
$$
We may choose here
$$
v(x) \ses x_1-x_2+i(x_3-x_4)\ess .
\eqno 5.6
$$
As in Section 3, we set $\om=e^{2\pi i/3}$ and
$$
y_1 = x_1x_2+x_3x_4
\ess\ess ,\ess\ess 
y_2 = x_1x_3+x_2x_4
\ess\ess ,\ess\ess 
y_3 = x_1x_4+x_2x_3\ess .
\eqno 5.7
$$
Moreover we let
$$
w_1 = x_1+x_2-x_3-x_4
\ess\ess ,\ess\ess 
w_2 = x_1+x_3-x_2-x_4
\ess\ess ,\ess\ess 
w_3 = x_1+x_4-x_2-x_3\ess .
\eqno 5.7
$$
This given, we find that in this case repetitive uses of Lagrange's Theorem 4.1 naturally leads us 
to the the following scheme:
$$
\matrix{
function& & expression & group\cr
\ess& \ess  & \ess & \ess\cr
\PP_o& \ses & (x_1-x_2)^2(x_1-x_3)^2\cdots (x_3-x_4)^2& S_4\cr  
\downarrow& \ess &\ess &2\downarrow\cr
\PP_1& \ses & (x_1-x_2)(x_1-x_3)\cdots (x_3-x_4)& A_4\cr  
\downarrow& \ess &\ess &3\downarrow\cr
\PP_2& \ses & y_1+\om y_2+\om^2 y_3& \{id,(1,2)(3,4),(1,3)(2,4),(1,4)(2,3)\}\cr  
\downarrow& \ess &\ess &2\downarrow\cr
\PP_3& \ses &w_2w_3+i(y_2-y_3)& \{id,(1,2)(3,4)\}\cr  
\downarrow& \ess &\ess &2\downarrow\cr
\PP_4& \ses & x_1-x_2+i(x_3-x_4)& \{id\}\cr}
\eqno 5.8
$$ 
The fact that $G_{\PP_o}= S_4$ is immediate since $\PP_o$ is a symmetric function of the roots
whose expression in terms of $e_1,e_2,e_2,e_4$ is given by the discriminant formula 2.8. Clearly
$\PP_1$ is invariant only under even permutations of the roots. So $G_{\PP_1}$ is simply the
alternating group $A_4$. To obtain $G_{\PP_2}$, we simply observe that, since $1,\om,\om^2$ are
distinct, $\PP_2$ is invariant only under those permutations of the roots that leave $y_1,y_2,y_3$
individually invariant. In other words $G_{\PP_2}=G_{y_1}\cap G_{y_3}\cap G_{y_3}$. 
This gives us the third entry in the fourth column of 5.8. To obtain $G_{\PP_3}$ we note that
for $y_2-y_3$ not to change we need each of $y_1,y_2,y_3$ to  remain unchanged, thus 
$G_{y_2-y_3}=G_{\PP_2}$. Now it is easily seen that
$$
G_{w_2w_3}\ses \{id,(1,2),(3,4),(1,2)(3,4)\}
$$
and thus we must have
$$
G_{\PP_3}\ses G_{w_2w_3}\cap G_{\PP_2}\ses  \{id,(1,2)(3,4)\}\ess 
$$
as asserted in 5.8.
\sas

We can easily see that 
$$
\PP_1^2\ses\PP_o\ess\ess {\rm and}\ess \ess \PP_4^2\ses \PP_3
\eqno 5.9
$$
We can painlessly check on the computer that 
$$
\PP_2^3\ses {\scriptstyle{3\over 2}}{\root 2\of 3}\ssp \PP_1\sps J\ess ,
\eqno 5.10
$$
where
$$
J\ses \sh(2y_1-y_2-y_3)(2y_2-y_1-y_3)(2y_3-y_1-y_2)\ess .
$$
Note next that we have the coset decomposition
$$
G_{\PP_2}\ses G_{\PP_2}\sps (1,3)(2,4) G_{\PP_2} 
$$
Thus the {\ita conjugate} of $\PP_3$ in $G_{\PP_2}$ is
$$
\PP_3'\ses (1,3)(2,4)\PP_3\ses -w_2w_3+i(y_2-y_3)\ess .
$$
We can now immediately conclude from Theorem 4.2 that the coefficients of the polynomial
$$
(t-\PP_3)(t-\PP_3')=\big(t-i(y_2-y_3)\big)^2- w_2^2w_3^2=
t^2- 4i(y_2-y_3) t - 4(y_2-y_3)^2- w_2^2w_3^2
\ses 
\eqno 5.11
$$
must be in $Rat[\CF;e_1,e_2,e_3,e_4,\PP_2]$. It is interesting to see what they actually turn out to be.
For instance we can write
$$
w_2^2w_3^2\ses {(w_1w_2w_3)^2\over w_1^2}\ess .
\eqno 5.12
$$
Now from 3.10 we get that
$$
w_1^2\ses e_1^2-4e_2+4y_1\ess .
\eqno 5.13
$$
Since $G_{\PP_2}\con G_{y_1}$, we know that $y_1$ should be in $Rat[\CF;e_1,e_2,e_3,e_4,\PP_2]$.
In fact, we can easily verify that
$$
y_1\ses {\scriptstyle {1\over 3}}\ssp (e_2+\PP_2+{\overline{\PP}}_2)\ess .
\eqno 5.14
$$  
In case we might worry that the complex conjugate ${\overline{\PP}}_2$ may not be in  $Rat[\CF;e_1,e_2,e_3,e_4,\PP_2]$,
I will quickly point out that we have the identity
$$
{\overline{\PP}}_2\ses {e_1^2- 3e_2\over \PP_2}\ess .
$$
Similarly, we can easily derive that
$$
i(y_2-y_3)\ses {1\over {\root 2\of 2}}\ssp (\PP_2- {\overline{\PP}}_2)\ess .
$$
Finally combining 5.11,5.12,5.13,5.14 we get that the function 
$$
\PP_3^*\ses \PP_3-i(y_2-y_3)
$$
satisfies the equation
$$
(\PP_3^*)^2\ses {w_1^2w_2^2w_3^2\over e_1^2-4e_2+4(e_2+\PP_2+{\overline{\PP}}_2)/3}\ess .
$$
Thus we have again obtained the desired $n!$-valued function by constructing a sequence of
functions $\PP_i(x)\in Rat[\CF;\xon]$ satisfying the recurrence
$$
\PP_i(x)^{p_i}\in Rat[\CF;\eon,\PP_{i-1}]
\bigsp (\ess {\rm with }\ess\ess \PP_o\in Rat[\CF;\eon]\ess ) 
\eqno 5.14
$$
Where, exponents $p_1,p_2,p_3\ldots$ give the prime factorization of $n!$.
Moreover we also have verified that the additional condition
$$
G_{\PP_i}\con G_{\PP_{i-1}} 
\eqno 5.15
$$
holds throughout.
\sas

This given, Lagrange was convinced (and, with hindsight, he was right),
that if the quintic should be solvable by radicals, then the $5!$-valued function
of its roots
must be reachable by a sequence of steps as given in 5.14 and 5.15.
However in trying to reproduce the same scheme for the quintic equation
he run into insurmountable difficulties. His final efforts were towards constructing
a rational function of the roots $x_1,x_2\,\ldots ,x_5$ that 
took less than $5$ values or one that took exactly five values but was the solution
of a binomial quintic. The reason for this search was his need to reduce
the construction of the roots to the solution of an equation of degree less than
five or to an equation of degree $5$ he could solve by taking radicals. Of 
course he was allowing the possibility,  of
having to solve an equation of degree greater than five as long as this
equation, as in the case of the cubic, could be reduced to the solution
of an equation of lesser degree for a power of one of its roots. However, this power
would then have to be a $k$-valued function for some $k<5$.
The best he could do was to produce a $6$-valued function. He concludes
his work (see []) by saying that although he had not tried every possibility,
the search (without {\small MAPLE} or {\small MATHEMATICA}) was considerably time consuming...
and he had no further time to spend in the search for something 

\centerline{\ita ``whose existence is very much in doubt''}

He had again the correct suspicion! The non existence of what he was looking for
was proved (30 years later) by Ruffini and generalized for $n>5$ by Cauchy.
\sas

\noindent{\bol Remark 5.1}

It is interesting in this connection to observe that there is a neat 
representation theoretical reason why there is no $k$-valued function
of the roots of $E_n(t)$ when $2<k<n$ and $n\geq 5$. In fact, 
the action of $S_n$ on the stabilizer of any $k$-valued function $\PP$ induces
a permutation representation with only one occurrence of the trivial.
If the remaining irreducible constituents are  all sign representations then $\PP$ is only $2$ valued.  
So for $k>2$ this representation would have to have a irreducible constituent of degree
$\geq 2$. But for $n\geq 5$ except for the trivial and the sign representation all
the other irreducible representations have dimension $\geq n-1$.
So $k\leq 2$ or $k\geq n$. At any rate we can also give this result an elementary proof. 
\sas

It is important to see at this point what are the implications of the existence of 
a sequence of functions $\PP_i$ satisfying 5.14 and 5.15. To this end we terminate
this section with a collection of results which shed considerable light on the
difficulties encountered by Lagrange in completing his program.
We should note that although some of the arguments that follow use the ``group''
structure of stabilizers, which is one of the main discoveries of Galois, most
of the calculations are actually due to Lagrange. This apparent paradox 
illustrates in a remarkable way how close Lagrange got to discover some of
the main points of ``Galois'' theory.

\noindent{\bol Theorem 5.1}

{\ita Let $\PP$ and $\PS$ be in $ Rat[\CF;\eon]$. Suppose that
$$
\eqalign{
&a)\ess\ess G_\PS\con  G_\PP\ess\ess but\ess\ess G_\PS\neq  G_\PP\ess\ess , and \cr
&b)\ess\ess \PS^p\in Rat[\CF;\eon,\PP]\cr }
\eqno 5.16
$$
where $p$ is prime. Then the decomposition of $G_\PP$ into left cosets
of $G_\PS$ may be written in the form
$$
G_\PP\ses G_\PS + \ggg G_\PS+\ggg^2 G_\PS+\cdots + \ggg^{p-1} G_\PS
\eqno 5.17
$$
with $\ggg$ a $p$-cycle that commutes with $G_\PS$. Moreover, the
conjugates of $\PS$ in $G_\PP$ can be written in the form
$$
\PS_i=\ggg^{i-1}\PS = \om^{i-1}\PS \bigsp (\ess with \ess\ess  \om=e^{2\pi i/p}\ess )
\eqno 5.18
$$
}

\sas

\noindent{\bol Proof}

Since by assumption $G_\PS\neq  G_\PP$ we have 
$$
G_\PP\ses \tau_1G_\PS+ \tau_2G_\PS+\cdots +\tau_kG_\PS\ess .
$$
with $k>1$. Set $\om=e^{2\pi i/p}$ and let
$$
Q(t)\ses \prod_{i=0}^{p-1}\ssp (t-\om^i \PS)\ess .
\eqno 5.19
$$
From 5.16 b) we deduce that
$$
Q(t)\ses t^p\sms R(x)
$$
with $ R(x)\in  Rat[\CF;\eon,\PP]$. In particular for any $\sig\in G_\PP$ we must have
$$
\sig Q(t)\ses \prod_{i=0}^{p-1}\ssp (t-\om^i \sig \PS)\ses Q(t)\ses \prod_{i=0}^{p-1}\ssp (t-\om^i \PS)
\ess .
$$
and this gives
$$
\sig\PS\ses \om^{i(\sig)}\ssp \PS
\eqno 5.20
$$
for some $0\leq i(\sig)\leq p-1$. This implies that the conjugates of $\PS$ 
are all multiples of $\PS$ by powers of $\om$. Thus
$$
k\leq p\ess .
$$  
Now $k>1$ gives that $\tau_2\PS\neq \PS$ so $\ess i(\tau_2)\neq 0\ess (mod \ess p)$. 
But since $p$ is
prime, the successive powers $(\om^{i(\tau_2)})^s$ for $s=0,1,\ldots ,p-1$ are simply $1,\om,\ldots
,\om^{p-1}$ in some order. 
Thus 
$$
\PS,\tau_2\PS,\tau_2^2\PS\ldots ,\tau_2^{p-1}\PS
$$
are all distinct. This gives 
$$
p\leq k\ess 
$$
Thus  $k$ must be equal to $p$ and
$$
\tau_2^p=id\ess .
$$
Moreover, we could have indexed our
coset representatives $\tau_i$ to be successive powers of an element $\ggg\in G_\PP$ 
for which $\ggg \PS=\om \PS$. This gives 5.17. With this choice we have
$$
\PS_i\ses \ggg^i\PS\ses \om^{i-1}\ssp \PS\bigsp (\ess i=1,2,\ldots ,p\ess ) 
$$
Finally, we see that $\ggg \PS=\om\PS$ gives
$$
\ggg^{-1} h\ggg\ssp  \PS\ses \PS \bigsp(\ess \forall\ess\ess h\in G_\PS)\ess .
$$
In other words
$$
\ggg^{-1} \ssp G_\PS \ggg \ses G_\PS
\eqno 5.21
$$
This completes our proof.
\sa

\noindent{\bol Remark 5.2}

In modern terminology this theorem simply says that 5.16 a) and b) imply that
$G_\PS$ is a normal subgroup $G_\PP$  and the quotient
$G_\PP/G_\PS$  is isomorphic to the group   of integers $mod\ssp p$.
To express these two properties  we shall here and after write
$$
G_\PS\nor_p G_\PP \ess .
\eqno 5.22
$$
Thus the possibility of solving the general equation $E_n(t)$ by radicals as was done
with the cubic and the quartic according to the scheme given by 5.14 and 5.15
requires that we should have a sequence of subgroups $G_i\con S_n$ (for $i=0,1,\ldots ,k)$
satisfying the conditions
$$
\{id\}= G_m \nor_{p_m} G_{m-1}\nor_{p_{m-1}} \cdots G_2\nor_{p_2} G_1\nor_{p_1} G_o= S_n
\eqno 5.23
$$
We shall soon see that this is impossible for $n\geq 5$. However, we must first
establish some basic facts about $S_n$ and its subgroups.  
\sap

\noindent{\bol Proposition 5.1}

{\ita Let $G$ be a subgroup of $S_n$ which contains all $3$-cycles and let $H$
be a subgroup $G$. Suppose that for some $\ggg\in G$ we have
$$
\eqalign{
&a)\ess\ess  G=H+\ggg H+\cdots +\ggg^{p-1}H\ess ,\cr
&b)\ess\ess  \ggg H \ses H\ggg\ess .\cr
}
\eqno 5.24
$$
Then for $n\geq 5$ also $H$ contains all $3$-cycles.}
\sa

\noindent{\bol Proof}

Condition a) in 5.24 gives that we can write
$$
(1,2,3)\ses \ggg^i a\ess\ess\ess ,\ess\ess\ess 
(3,4,5)\ses \ggg^j b
$$
for some $0\leq i,j<p$ and $a,b\in H$. On the other hand b) in 5.24 gives
that we can write
$$
(1,2,3)(3,4,5)\ses \ggg^i a\ggg^j b=\ggg^{i+j} \ssp c
\ess\ess\ess ,\ess\ess\ess 
(3,2,1)(5,4,3)\ses a^{-1}\ggg^{-i}b^{-1}\ggg^{-j}= c'\ggg^{-i-j}
$$
for some $c,c'\in H$. Thus 
$$
(2,5,3)\ses (3,2,1)(5,4,3) (1,2,3)(3,4,5)\ses c'\ggg^{-i-j}\ggg^{i+j} \ssp c\in H\ess .
$$
Since the same manipulations can be carried out when $1,2,3,4,5$ are replaced by
any $5$ distinct integers $i_1,i_2,i_3,i_4,i_5$, we see that $H$ must
contain all $3$ cycles as asserted.
\sa  

\noindent{\bol Proposition 5.2}

{\ita If a subgroup $G\con S_n$ contains all $3$-cycles then $G=S_n$ or  $G$
is the alternating group $A_n$}
\sa

\noindent{\bol Proof}

Note that we have
$$
\eqalign{
(1,2)(2,3)&\ses (1,2,3)\cr
(1,2)(3,4)&\ses (1,2)(2,3)(2,3)(2,4)\ses (1,3,2)(2,4,3) \cr}
$$
Thus the product of two $2$-cycles can always be expressed as either a $3$-cycle or as
a product of two $3$-cycles. This implies that every even permutation is a product of
$3$-cycles. Thus under our hypothesis $G\supseteq A_n$. Now if $G$ contains a single 
permutation $\sig$ not in $A_n$ (that is an odd permutation), the identity
$$
S_n\ses A_n\sps \sig A_n
$$
yields that $G=S_n$ as desired.
\sa

Combining these two propositions we derive that we cannot have 5.23 when $n\geq 5$.
In fact the string of inclusions in 5.23 must stop with $G_1=A_n$ and can go no further!
We can thus state 
\sa

\noindent{\bol Theorem 5.2} ( Ruffini)

{\ita For $n\geq 5$ the general equation $E_n(t)=0$ cannot be solved by radicals
by successive adjunctions of rational functions $\PP_i\in Rat[\CF,\xon]$ satisfying the
recursions in 5.14 and 5.15 as was possible for the cubic and the quartic.} 
\sa

We should note that this result doesn't quite put to rest the possibility of solving
the general equation by adjoining ``radicals''. This is because
there are still some unjustified hypotheses in the present setting.
First of all, in our original definition (see 4.15) of {\ita solution by radicals} 
we did not have the extra condition 5.15. As we shall see this is a relatively
minor point, easily fixed in the Galois setting. Considerably more troubling is that 
in this definition we require that each newly constructed $\PP_i$ should turn out to 
be a {\ita rational function in the roots}. What if we allow such steps as taking a $p^{th}$
root of $x_1x_2\cdots x_n$? Can we then construct a solution? These doubts would
be removed if we could show that the existence of a sequence of steps as in
5.14, 5.15 without the further assumption that $\PP_i\in Rat[\CF,\xon]$ implies the existence
of the same sequence of steps  with this assumption {\bol satisfied}.
This is in essence what Abel succeeded in showing. We shall not present Abel's arguments
here since what must be proved to remove the remaining doubts can be done in an easier 
way in the Galois setting. We shall thus terminate our treatment of Lagrange's
``Galois'' Theory
with one final result which in some sense reverses Theorem 5.1.  
\sa

\noindent{\bol Theorem 5.3}

{\ita If $\PS,\PP\in Rat[\CF;\xon]$ and 
$$
G_\PP\ses G_\PS + \ggg G_\PS+\ggg^2 G_\PS+\cdots + \ggg^{p-1} G_\PS
\eqno 5.25
$$
with $\gg^p=id$ and
$$
\ggg G_\PS\ses G_\PS\ggg\ess . 
\eqno 5.26
$$
Then we can find $\TT\in  Rat[\CF;\xon]$ such that}
$$
\eqalign{
&1)\ess\ess G_{\TT}\son G_\PS \ess\ess  ( \TT\in Rat[\CF;\eon,\PS]) ,\cr
&2)\ess\ess \TT^p \in Rat[\CF;\eon,\PP].\cr}
\eqno 5.27
$$
{\bol Proof}

Set
$$
\TT\ses \PS  +\om\ggg \PS +(\om\ggg)^2\PS +\cdots +(\om\ggg)^{p-1}\PS \ess .\bigsp 
(\ess \om=e^{2\pi i/3}\ess) 
\eqno 5.28
$$
It is easily verified that 
$$
\ggg\TT\ses \om^{-1}\TT\ess .
\eqno 5.29
$$
Moreover, 5.26 gives that each of the conjugates $\PS_i=\ggg^{i-1}\PS$
has the same stabilizer as $\PS$. Thus from the definition 5.28 we get that
$$
G_\TT\son G_\PS
\eqno 5.30
$$
and therefore (again by Theorem 4.1) we must have
$$
\TT\in Rat[\CF;\eon,\PS]
$$
Finally, 5.30 together with 5.25 and 5.26 implies that the polynomial
$$
Q(t)\ses \prod_{i=1}^p (t-\ggg^{i-1}\TT)\ses \prod_{i=1}^p(t-\TT/\om^{i-1})\ses 
t^p-\TT^p
\eqno 5.31
$$
is left invariant by every element of $G_\PP$. Thus from Theorem 4.1 
(with $\PS\rightarrow \PP$ and $\PP \rightarrow \TT$) we derive that
$$
\TT^p\in Rat[\CF;\eon,\PP]\ess ,
$$
as desired. 
\sa

\noindent{\bol Remark 5.3}

We should point out that the inequality in 5.30 cannot in general be sharpened to an equality.
Indeed, in our construction  of the roots of the quartic, the function
$$
\PP_3=w_2w_3+i(y_2-y_3)
$$
satisfies the conditions in 5.25 and 5.26 with $\PS\rightarrow\PP_3$,
$\PP\rightarrow\PP_2$ $p=2$. In this
case we have
$$
\ggg=(1,3)(2,4)\ess\ess ,\ess\ess \om=-1
$$
and 
$$
\TT= \PP_3+\om \ggg\PP_3=2 w_2w_3\ess .
$$
As we have seen
$$
G_{w_2w_3}\ses \{id,(1,2),(3,4),(1,2)(3,4)\}
$$
while
$$
G_{\PP_3}\ses \{id,(1,2)(3,4)\}\ess . 
$$
\sa

Nevertheless, in a typical application of Theorem 5.3 we may want to obtain $\TT$
by extracting a $p^{th}$ root of a function in $Rat[\CF;\eon,\PP]$ then 
recover $\PS$ by rational operations involving $\TT$ and possibly other known functions. 
Now this is provided by the following beautiful identity. 
\sa

\noindent{\bol Theorem 5.4} (Lagrange)

{\ita If $\PS,\PP\in Rat[\CF;\xon]$ and we have 5.25, 5.26 with $\ggg^p=id$
and $p$ a prime. Then
$$
\PS\ses (c_o+c_1\TT+\cdots +c_{p-1}\TT^{p-1})/p
\eqno 5.32
$$
where all the coefficients $c_i$ are in $Rat[\CF;\eon,\PP]$}.
\sas

\noindent{\bol Proof}

For convenience set $\PS_s=\ggg^s\PS$
$$
\TT_i\ses \sum_{s=0}^{p-1}\ssp \om^{si} \PS_s\ess .
$$
Then it is easily verified that
$$
(1-\om^i\ggg)\TT_i\ses 0\ess .
\eqno 5.33
$$
Moreover, since for $p$ a prime and any $0<s<p$ we have 
$$
1 \sps \om^s\sps \om^{2s}+\cdots +\om^{(p-1)s}\ses 0\ess ,
$$
we immediately derive that
$$
\PS\ses (\TT_1+\TT_2+\cdots + \TT_p)/p\ess .
\eqno 5.34
$$
On the other hand, 5.33 gives that the ratios
$$
c_i\ses \TT_i/\TT^i
\eqno 5.35
$$
are all invariant under $\ggg$. Since as we have seen 5.26 gives that all
the conjugates $\PS_i$ of $\PS$ in $G_\PP$ have $G_\PS$ as stabilizer, we deduce
that these ratios are stable under $G_\PP$ and (again from Thorem 4.1) we must
conclude that they are all in $Rat[\CF;\eon,\PP]$. This given, we see that by
substituting 5.35 in 5.34 we get 5.32 with the desired properties.
\vfill\supereject
\def \om {\omega}
\def\CF {{\cal F}}
\def\CE {{\cal E}}
\def\CA {{\cal A}}
\def\CB {{\cal B}}

\def\RFX {Rat[{\cal F};x_1,\ldots ,x_n]}
\def\PFX {{\cal F}[x_1,\ldots ,x_n]}
\def \xsig {x_{\sig_1},x_{\sig_2},\ldots ,x_{\sig_n}}
\def \PP {\Phi}
\def \PS {\Psi}
\def \SYM {Sym[\CF;\xon]}
\def \RSYM {Ratsym[\CF;\xon]}
\def \part {\vdash}
\def \DD {\Delta}
\def \omi {\hskip -.01truein\ ^{(i)}}
\def \scos {\ssp :\ssp}
\def \bup {\hskip -.01truein}
\def \om {\omega}
\def \eon {e_1,e_2,\ldots ,e_n}
\def \SQR {{\root 2 \of R}}
\def \sh {{\scriptstyle{1\over 2}}}
\def \shf {{\scriptstyle{1\over 4}}}
\def \AO {{\root 3 \of {-\sh q+ \SQR}}}
\def \BO {{\root 3 \of {-\sh q- \SQR}}}
\def \aon {\aaa_1, \aaa_2,\ldots ,\aaa_n}

\def \TPP {{\tilde \PP}} 
\def \TE {{\tilde E}} 
\def \TT {\Theta} 
\def \TP {{\tilde P}} 
\def \TQ {{\tilde Q}}
\def \tv  {{\tilde v}}
\def \TPS  {{\tilde \PS}}
\def \RA {\ssp \Longrightarrow\ssp}

\def \nor {\triangleleft}
\def \son {\supseteq}
\def \res {\TP_1(t)}
\def \fora {\bigsp (\ess \forall \ess\ess \ggg\in G\ess )}
\sa
\noindent{\bol 6. Galoisian ``Galois Theory.}
\sa

Lagrange's pursuits were brought to a conclusion by Abel around 1829 when Galois started his
investigations. Galois' starting point may have been the idea that although there may be
no general formulas for the roots of $E_n(t)$ for $n\geq 5$ which involved only extraction of roots,
it may still be possible to find them for some special equations. This possibility had already
emerged in the work of Gauss and Abel but it is not clear to what extent Galois had been aware of 
their work.   We may say that his discoveries stemmed from a natural extension of Lagrange's original methods. 
To present Galois'
results, we shall need some additional notation and definitions. \sa

Our basic ingredients here will be two fields $\CF$ and $\CE$, with $\CF$ a proper subfield of $\CE$. 
The equation to be solved will be written as
$$
\TE_n(t)\ses (t-\aaa_1)(t-\aaa_2)\cdots (t-\aaa_n)\ses 0\ess .
\eqno 6.1
$$
where $\aon$ are distinct and, unless explicitely mentioned, will remain unchanged
throughout our presentation. We shall also assume hereafter that
$$
\eqalign{
&1)\ess\ess \aon\ess \in \ess \CE \ess\ess {\rm and }\cr
&2)\ess\ess e_1(\aon ),e_2(\aon ),\ldots ,e_n(\aon )\ess \in \ess \CF\cr}
\eqno 6.2
$$
We see that we shall have to work here with functions $\PP(\xon)$
of the independent variables $\xon$ and at the same study their values when 
when $\xon$ are replaced by $\aon$. As clarity requires, these values will 
be represented by any of the symbols below
$$
\PP(\aon )\ses \PP(\aaa ) \ses \TPP\ess .
\eqno 6.3
$$
For a given $\sig\in S_n$ we shall also use the symbols
$$
\sig\PP(\aon )\ses \sig\PP(\aaa ) \ses \sig\TPP\ess 
$$
to denote the value $\PP(\aaa_{\sig_1},\aaa_{\sig_2},\ldots ,\aaa_{\sig_n})$.
\sas

This given, one of the fundamental differences between the Galoisian and Lagrangian setups 
is that although a given $\PP\in Rat[\CF;\xon]$ and all its images $\sig\PP$ are well defined
as elements of $ Rat[\CF;\xon]$ some of the  values $\sig\TPP$ may make no sense at at all.
An example in point is the rational function
$$
\PP(x)\ses {1\over x_1^2+x_2^2-x_3^2}
$$
when $n=3$ and $\aaa_2^2+\aaa_3^2-\aaa_1^2=0$.
\sas

In summary in the Galois setting, we have to be careful with denominators! We shall avoid the
problem by dealing hereafter only with polynomials functions of the roots. As we shall see this is not a serious 
restriction, and with it, most of the results of Lagrange Theory can be extended to the Galois
setting with nearly identical proofs.
\sas

Another important difference is that for some $\PP\in \PFX$ there may be more permutations 
of the roots of $\TE_n(t)$ that leave  $\PP(\aon )$ unchanged than there are in $G_\PP$.
We shall take account of this difference by setting for $\PP\in \PFX$
$$
S_\PP\ses \{\sig\in S_n\ess :\ess \PP(\aaa_{\sig_1},\aaa_{\sig_2},\ldots ,\aaa_{\sig_n})=
\PP(\aaa_1,\aaa_2,\ldots ,\aaa_n)\ess \}
\eqno 6.4
$$
Note that $S_\PP$ should be considered a property of $\PP$ and not a property of the value $\PP(\aaa )$.
Moreover we should point out that in general this collection of permutations may not even be a group!
For instance for the equation
$$
\TE_3(t)\ses \ses (t+1)(t+i)(t-i)\ses t^3+t^2+t+1\ess ,
$$
we can take $\CF$ and $\CE$ to be  the fields of rational and  complex numbers respectively. Now,
if we label the roots by setting 
$$
-1=\aaa_1
\ess\ess ,\ess \ess
i=\aaa_2
\ess\ess ,\ess \ess
-i=\aaa_3 
$$
then for $\PP(x_1,x_2,x_3)\ses x_2^2$
we have (in cycle notation)
$$
G_\PP\ses \{id\scs (1,3)\}\ess .
$$
On the other hand the permutations of $\aaa_1,\aaa_2,\aaa_3$ that leave $\TPP$ unchanged form the set
$$
\{id\scs (1,3)\scs (2,3)\scs (1,2,3)\}\ess 
$$
This is not a group since $(1,2)=(1,2,3)(2,3)$, yet 
$$\PP(\aaa_1 ,\aaa_2 ,\aaa_3)=\aaa_2=-1 \ess\ess {\rm and} \ess\ess \PP(\aaa_2,\aaa_1,\aaa_3)=\aaa_1^2=1\ess .
$$
To take account of this possibility we shall say that a given $\PP\in \PFX$ is  {\ita Galois} if and only if 
$$
S_\PP \ses  G_\PP
$$

We are now in a position to proceed with our treatment.
\sa

\sa

\noindent{\bol Theorem 6.1}

{\ita For every subgroup $ G\con S_n$ we have a Galois $\PP\in \CF[\xon]$ such that}
$$
S_\PP=  G\ess .
\eqno 6.5
$$
{\bol Proof}

We start with the case $G=S_n$. Here can take again a linear function 
$$
v(\xon)\ses m_1 x_1+m_2 x_2+\cdots + m_n x_n
\eqno 6.6
$$
as in 4.16, but we must be a bit more careful in choosing the coefficients $m_i$. To this end note that
if $m_1,m_2,\ldots , m_n$ are chosen to be integers in the interval $[0, M]$, then
for any given pair of distinct permutations $\sig,\tau\in S_n$ the equation
$$
m_1 \aaa_{\sig_1}+m_2 \aaa_{\sig_2}+\cdots + m_n \aaa_{\sig_n}\ses 
m_1 \aaa_{\tau_1}+m_2 \aaa_{\tau_2}+\cdots + m_n \aaa_{\tau_n}
$$
can have at most $(M+1)^{n-1}$ distinct solution vectors $(m_1,m_2,\ldots ,m_n)$.
Thus to assure that 
$$
\sig v(\aaa) \neq \tau v(\aaa) \ess\ess\ess (\ess \forall \ess\ess i\neq j\ess )
$$
we need to avoid at most ${n! \choose 2}\times (M+1)^{n-1}$ vectors.  
However, when $M+1>{n!\choose 2}$ there will remain some  $n^{tuples}$ for us to choose and satisfy our
requirement that $S_v=S_n$. Having made one such choice of $(m_1,m_2,\ldots ,m_n)$, the desired $\PP$ for any given
subgroup $G$ can be readily produced. In fact, we can show that we can set
$$
\PP\ses \prod_{\sig\in G}\ssp (N-\sig v(x))
\eqno 6.7
$$
where $N$ is a suitably chosen integer. To see this note first that the form of 6.7 guarantees
that whatever $N$ we choose we shall have at least $ S_\PP\supseteq G$. Now let
$$
S_n\ses \tau_1G \sps \tau_2G\sps \cdots \sps \tau_kG
$$   
be the decompostion of $S_n$ into left cosets of $G$. This given, our choice of $m_1,m_2,\ldots, m_n$ assures
that the polynomials 
$$                
\TP_i(t)\ses \tau_i\prod_{\sig\in G}\ssp (t-\sig v(\aaa))\ses \prod_{\sig\in G}\ssp (t-\tau_i\sig v(\aaa))
$$
have no roots in common. Since they all have degree $|G|$, the equation
$$
\TP_i(t)=\TP_j( t)
$$
for $i\neq j$ can then have at most $|G|$ solutions. Thus if we want an integer $N$ which
gives $\TP_i(N)\neq \TP_j(N)$ for all $i\neq j$ we need avoid at most ${k\choose 2}\times |G|$
values. Clearly we can find such an $N$ in the interval $[0,M]$ as soon as $M>{k\choose 2}\times |G|$.
This completes our argument.
\sa

\noindent{\bol Remark 6.1}

All the constructions, proofs and definitions in this section will use an $n!$-valued
Galois function
$$
v(x)\ses m_1x_1+m_2x_2+\cdots +m_nx_n
$$
which must remain unchanged throughout the rest of the section. 
We must therefore make sure that some of the objects we introduce, such as for instance
the ``Galois Group'' of our equation $\TE_n(t)$ do not depend on the choice of
$m_1,m_2,\ldots ,m_n$. This is one of the prices we have to pay for not following
the abstract approach. However, we believe that this will be well compensated by 
the additional insights that our insistence on explicit constructions will provide.
\sa

We begin by showing that the values of every polynomial $\PP(x)$ are in fact  
polynomials in the values of $v(x)$. More precisely we have
\sas

\noindent{\bol Proposition 6.1}

{\ita For any $\PP\in \PFX$ we can construct a polynomial $\TT_\PP(t)\in \CF[t]$
such that for any $\ggg\in S_n$ we have}
$$
\ggg\TPP\ses \TT_\PP\big(\ggg \tv\big)\ess .
\eqno 6.8
$$    
\noindent{\bol Proof}

As in section 4  (see 4.19 and 4.20) we let
$$
Q(t)\ses \sum_{\sig\in S_n}\ssp \sig \PP\prod_{\tau\in S_n}\bup ^{(\sig)}\ssp (t-\tau v)
\ess\ess ,\ess\ess
P(t)\ses \prod_{\tau\in S_n} (t-\tau v)\ess .
\eqno 6.9
$$
Since both $Q$ and $P$ are by construction $S_n$-invariant, their coefficients are polynomials
in $\CF[\eon]$. The hypothesis in 6.3 2) then yields that the polynomials 
$$
\TQ(t)\scs \TP(t)\scs \TP'(t)
$$ 
have coefficients in $\CF$. Setting $t=\ggg \tv$ in
$$
\TQ(t)\ses \sum_{\sig\in S_n}\ssp \sig \TPP\prod_{\tau\in S_n}\bup ^{(\sig)}\ssp (t-\tau \tv)
$$
gives
$$
\TQ(\ggg \tv )\ses \ggg \TPP\ssp \TP'(\ggg \tv)\ess .
\eqno 6.10
$$
Note further that since by construction $v$ takes $n!$ distinct values,
the polynomials $\TP$ and $\TP'$ have no common root. We can thus apply the Berlekamp algorithm
and construct a pair of polynomials $p(t),q(t)\in \CF[t]$ such that
$$
p(t)\TP'(t)\sps q(t)\TP(t)\ses 1
$$
Setting $t=\ggg\tv$ in this equation yields
$$
p(\ggg\tv)\TP'(\ggg\tv) \ses 1\ess .
$$
Multiplying both sides of 6.10 by $p(\ggg\tv)$ and using this equation we finally get
$$
p(\ggg\tv)\TQ(\ggg \tv )\ses \ggg \TPP\ess ,
$$
and this gives 6.8 with
$$
\TT_\PP (t)\ses p(t)\TQ(t )\ess .
$$
\sa

It develops that the polynomial 
$$
\TP(t)\ses \prod_{\tau\in S_n}\ssp (t-\tau\tv)
\eqno 6.11
$$ 
plays a crucial role  in our development. It may or may not be reducible in $\CF[t]$. If it
is, we can write its factorization into irreducibles in the form
$$
\TP(t)\ses \prod_{\tau\in T_1}(t-\tau\tv)\ess \prod_{\tau\in T_2}(t-\tau\tv)\cdots 
\prod_{\tau\in T_m}(t-\tau\tv)\ses \TP_1(t)\TP_2(t)\cdots \TP_m(t)\ess ,
$$
where $T_1,T_2,\cdots ,T_m$ are disjoint subsets, and $T_1$ is the subset that contains the identity permutation. 
Now we have the following crucial fact.
\sa

\noindent{\bol Proposition 6.2}

{\ita The subset $T_1$ is a group}
\sas

\noindent{\bol Proof}

Since by definition $T_1$ contains the identity, we only need to show that if $\eta , \xi\in T_1$
then their product $\eta \xi $ is also in $T_1$. To this end, we use Proposition 6.1 and construct
the polynomial $\TT(t)\in \CF[t]$ that gives 6.8 for $\PP=\xi v$. We shall thus have
$$
\ggg \xi \tv\ses \TT(\ggg \tv )\bigsp (\ess \forall \ess\ess \ggg\in S_n)\ess .
\eqno 6.12
$$
Now note that by hypothesis
$$
\TP_1(\xi\tv) \ses \prod_{\tau\in T_1}(\xi\tv-\tau\tv)\ses 0.
$$
In particular, using 6.12 for $\ggg=id$ we can rewrite this in the form
$$
\TP_1\big(\TT(\tv)\big)\ses 0\ess .
$$
But this says that the polynomials $\TP_1(t)$ and $\TP_1\big(\TT(t)\big)$ have a root in common.
Since they are both in $\CF[t]$ and $\TP_1(t)$ is irreducible in $\CF[t]$,
we deduce that $\TP_1\big(\TT(t)\big)$
must vanish for all the other roots of $\TP_1(t)$. In particular we must have
$$
\TP_1\big(\TT(\eta \tv)\big)\ses 0
$$
Now using 6.12 with $\ggg=\eta$ this may yet be rewritten as
$$
\TP_1(\eta\xi\tv)\ses 0\ess .
$$
But this implies that the permutation $\eta\xi$  lies  in $T_1$ as well. 

\hfill Q.E.D.
\sa
 
We shall hereafter denote $T_1$ by $G$ and refer to it as the {\ita Galois Group} of 
$\TE_n(t)$. The polynomial $\TP_1(t)$ itself will be referred to as a {\ita Galois resolvent}
of $\TE_n(t)$. For instance when $n=3$ and
$$
\TE_3(t)\ses t^3+t^2+t+1\ess ,
$$
we may take $\CF$ to be the field of rational numbers and $v=x_2-x_1$. This gives
$$
P(t)\ses \bigl(t^2-(x_2-x_1)^2\bigr)\bigl(t^2-(x_2-x_3)^2\bigr)\bigl(t^2-(x_3-x_1)^2\bigr)\ess .
$$
Now this can be rewritten as
$$
P(t)\ses t^6+(6e_2-2e_1^2)+(e_1^2-3e2)^2t^2-(x_1-x_2)^2(x_1-x_3)^2(x_2-x_3)^2
$$
and substituting $e_1=1$, $e_2=1$ and $e_3=-1$ we get (using formula 2.7)
$$
\TP(t)\ses t^6+4t^4+4t^2+16\ess . 
$$
Its irreducible factorization is
$$
\TP(t)\ses (t^2-2t+2)(t^2+2t+2)(t^2+4)\ess .
$$
Now the roots of $t^3+t^2+t+1$ are $-1,i,-i$ so if we label them $\aaa_1,\aaa_2,\aaa_3$
respectively, then the Galois resolvent is
$$
\TP_1(t)\ses t^2-2t+2\ses \bigl(t -(x_2-x_1)\bigr)\bigl(t -(x_3-x_1)\bigr)\ess ,
$$
and the Galois group reduces to
$$
G=\{ id, (2,3)\}\ess .
\eqno 6.13
$$ 
Note that if we had chosen $v=x_2-x_3$ then the Galois resolvent would have been
$$
t^2+4\ses \bigl(t-(x_2-x_3)\bigr)\bigl(t+(x_2-x_3)\bigr)
$$
and the Galois group would still be as in 6.13. 
It is easy to see from this example that as a subgroup of $S_3$, $G$ does depend on our labeling 
of the roots. Nevertheless we are going to show that as a group of
permutations of the set $-1,i,-i$, $G$ only depends on the equation $\TE_n(t)$ and the given
field $\CF$.
\sa

To this end, we need to introduce two classes of subgroups of $S_n$. We shall set
$$
\CA\ses \{\ssp H\con S_n \ssp :\ssp \PP\in \PFX\ssp \&\ssp S_\PP\son H \RA \TPP\in \CF\ssp \}
\eqno 6.14
$$
and
$$
\CB\ses \{\ssp H\con S_n \ssp :\ssp \PP\in \PFX\ssp \&\ssp \TPP\in \CF \RA  S_\PP\son H\ssp  \}\ess .
\eqno 6.15
$$
In words, a subgroup $H$ of $S_n$ belongs to $\CA$ if and only if any polynomial $\PP\in \PFX$
whose value is invariant under $H$ has its value $\TPP$ in $\CF$. In the same vein we can
say that a subgroup $H$ of $S_n$ belongs to $\CB$ if and only if any polynomial $\PP\in \PFX$
with its value $\TPP$ in $\CF$ must remain invariant (by value) under all elements of $H$.
It is immediate from the definitions 6.14 and 6.15 that for any two groups $H,K\con S_n$ we
have
$$
H\in \CA\ssp \&\ssp K\son H  \RA K\in \CA 
\ess\ess\ess\ess  {\rm and}\ess\ess\ess\ess  
H\in \CB\ssp \&\ssp K\con H  \RA K\in \CB 
\eqno 6.17
$$
In words, $\CA$ and $\CB$ are respectively   upper and lower ideals of subgroups of $S_n$
(under containement). Remarkably, we have the following basic fact
\sa

\noindent{\bol Theorem 6.2}

{\ita $\CA$ and $\CB$ are both principal ideals with $G$ as their unique extremal element.
That is }
$$
\CA\ssp \cap \ssp \CB \ses \{ G\}
\eqno 6.18
$$
\noindent{\bol Proof}

We start by proving that
$$
G\in \CA\ssp \cap \ssp \CB\ess .
\eqno 6.19
$$
Given a $\PP\in \PFX$ and using Proposition 6.1 we may write
$$
\ggg \TPP\ses \TT_\PP(\ggg \tv)\bigsp (\ess \forall \ggg\in S_n\ess )
\eqno 6.20
$$
Thus if $S_\PP\son G$ we have
$$
\TPP\ses {1\over |G|}\ssp \sum_{\ggg\in G}\ssp \TT_\PP(\ggg \tv)
$$
Now the right hand side of this expression is a 
symmetric polynomial  \footnote {(*)} {with coefficients in $\CF$}
in the roots of $\TP_1(t)$. This shows that $\TPP$ may equally be  expressed as a polynomial (*) 
in the coefficients of $\TP_1(t)$ which themselves are in $\CF$. This implies that $\TPP\in \CF$
and that 
$$
G\in \CA\ess .
$$
Conversely, if $\TPP=a\in \CF$ then using the same polynomial $\TT_\PP$ we may rewrite this as
$$
\TT_\PP(\tv)-a\ses 0\ess .
$$
But this says that the  polynomial
$$
\TT_\PP(t)-a\in\CF[t]
$$
has a root in common with $\res$. Thus it must vanish at all the other roots of $\res$.
That is we must have
$$
\ggg\TPP\ses \TT_\PP(\ggg\tv)\ses a\ses \TPP\bigsp (\ess \forall \ess\ess \ggg\in G\ess )\ess .
$$
This implies that $S_\PP\son G$ and that
$$
G\in \CB
\eqno 6.21
$$
To complete our argument we must show that $G$ is contained in all the other elements of $\CA$
and that $G$ contains all the other elements of $\CB$. 

Note that if $G_1\in \CA$, the polynomial
$$
Q(t)\ses \prod_{\ggg\in G_1}\ssp (t-\ggg\tv)\ess .
$$
(whose  coefficients are necessarily invariant under $G_1$) must belong to $\CF[t]$.
Since it has the root $\tv$ in common with $\res$ and $\res$ is irreducible $Q(t)$ must
be divisible by $\res$. This gives that
$$
G\con G_1\ess .
$$
Conversely, let $G_1\in \CB$. Consider the polynomial
$$
\PP(x)\ses \TP_1 \bigl(v(x)\bigr)\in \PFX\ess .
$$
Since its value $\TPP=\TP_1 (\tv)=0$ is clearly in $\CF$ it must remain invariant under
every element of $G_1$. That is we must have
$$
\TP_1 (\ggg\tv)=0
$$
for all $\ggg\in G_1$. This shows that
$$
G_1\con G
$$
and completes our proof.
\sa
\def \om {\omega}
\def\CF {{\cal F}}
\def\CE {{\cal E}}
\def\CA {{\cal A}}
\def\CB {{\cal B}}

\def\RFX {Rat[{\cal F};x_1,\ldots ,x_n]}
\def\PFX {{\cal F}[x_1,\ldots ,x_n]}
\def \xsig {x_{\sig_1},x_{\sig_2},\ldots ,x_{\sig_n}}
\def \PP {\Phi}
\def \PS {\Psi}
\def \SYM {Sym[\CF;\xon]}
\def \RSYM {Ratsym[\CF;\xon]}
\def \part {\vdash}
\def \DD {\Delta}
\def \omi {\hskip -.01truein\ ^{(i)}}
\def \scos {\ssp :\ssp}
\def \bup {\hskip -.01truein}
\def \om {\omega}
\def \eon {e_1,e_2,\ldots ,e_n}
\def \SQR {{\root 2 \of R}}
\def \sh {{\scriptstyle{1\over 2}}}
\def \shf {{\scriptstyle{1\over 4}}}
\def \AO {{\root 3 \of {-\sh q+ \SQR}}}
\def \BO {{\root 3 \of {-\sh q- \SQR}}}
\def \aon {\aaa_1, \aaa_2,\ldots ,\aaa_n}

\def \TPP {{\tilde \PP}} 
\def \TE {{\tilde E}} 
\def \TT {\Theta} 
\def \TP {{\tilde P}} 
\def \TQ {{\tilde Q}}
\def \tv  {{\tilde v}}
\def \TPS  {{\tilde \PS}}
\def \TG {{\tilde G}}
\def \RA {\ssp \Longrightarrow\ssp}

\def \nor {\triangleleft}
\def \son {\supseteq}
\def \res {\TP_1(t)}
\def \fora {\bigsp (\ess \forall \ess\ess \ggg\in G\ess )}

\noindent{\bol Remark 6.2}

Since the definitions 6.14 and 6.15 of the classes $\CA$ and $\CB$ only involve the given field $\CF$
and the roots of the equation $\TE_n(t)$, we see that one of the consequences of Theorem 6.2 is that also
$G$ itself only depends on $\CF$ and $\TE_n(t)$.  When in  our developments we keep the given equation fixed
and only vary the field, for simplicity, we shall use the notation $\CA(\CF)$, $\CB(\CF)$, $G(\CF)$  
and leave the dependence on ${\TE_n(t)}$ implicit. In all other cases we will indicate this
dependence with a subscript. We should also keep in mind that one of the immediate consequence
of our definition of a Galois group is that if $\CF\con\CF_1$ are two fields then we necessarily 
must have $G_{\TE_n(t)}(\CF_1)\con G_{\TE_n(t)}(\CF)$. 
\sa

In trying to extend Lagrange's Theorems 4.1 $\&$ 4.2 to the Galois setting we should be tempted to let the
Galois group $G$ play the role of $S_n$ in the arguments. However, given a polynomial $\PS\in \CF[\xon]$  
we may not be in a position to write down the coset decomposition in 4.5 with $S_n$ replaced by $G$
for the simple reason that we may not have $G_\PS\con G$. Using $G_\PS\cap G$ in place of $G_\PS$
doesn't get us anywhere for the simple reason that the intersection $G_\PS\cap G$ may only consist of
the identity permutation. Nor we can use  $S_\PS$ instead of $G_\PS$ for it may be too big
and as we have seen it may not even be a group. It develops that the optimal choice turns out to be the intersection
$ S_\PS\cap G$ which here and after will be denoted by $\TG_\PS$ and referred to as the {\ita Galois stabilizer} of 
$\PS$. In fact, $\TG_\PS$ is neither too small nor to big and remarkably it can be easily shown that
\sa

\noindent{\bol Theorem 6.3}

{\ita For any polynomial $\PS\in \CF[\xon]$ the  Galois stabilizer 
$$
\TG_\PS\ses S_\PS\cap G
$$
is a group}

\sas

\noindent{\bol Proof}

Since by construction $\TG_\PS$ contains the identity, we need only show that 
$$
\eta, \xi\in \TG_\PS\RA \eta^{-1}\xi \in \TG_\PS\ess .
\eqno 6.22 
$$
However, this is immediate.  In fact,  $\eta , \xi\in \TG_\PS$ implies that
$$
\xi  \TPS - \eta  \TPS= 0\ssp \in \CF\ess .
$$
But since $G\in \CB$ we must necessarily have
$$
S_{\xi\PS-\eta\PS} \son G
$$
In other words
$$
\ggg\xi  \TPS - \ggg\eta  \TPS\ses 0\ess .
\bigsp (\ess \forall \ess\ess \ggg\in G\ess)
$$
In particular this must hold true for $\ggg=\eta^{-1}$. That is 
$$
\eta^{-1}\xi\TPS\ses \TPS \ess ,
$$
which gives 6.22 as desired.
\sa

We are now in a position to state and prove the Galois version of Theorems 4.1 and 4.2.
\sas

\noindent{\bol Theorem 6.4}

{\ita If $\PS\&\PP\in \CF[\xon]$ and
$$
\TG_\PS\con \TG_\PP
\eqno 6.23 
$$
then  we can construct a polynomial $\TT(t)\in \CF[t]$ which gives}
$$
\ggg\TPP\ses \TT(\ggg\TPS) \fora
\eqno 6.24
$$
\sas

\noindent{\bol Proof}

Let
$$
G\ses 
\tau_1\ssp \TG_\PS\sps \tau_2\ssp \TG_\PS\sps\cdots\sps \tau_k\ssp \TG_\PS \bigsp (\ess \tau_1=identity\ess )
\eqno  6.25
$$
and set
$$
A(t)\ses  \sum_{i=1}^k\ssp \tau_i \TPP \ssp \prod_{j=1}^k\omi\ssp (t-\tau_j\TPS)
\ess\ess \scs \ess\ess
B(t)\ses \prod_{i=1}^k \ssp (t-\tau_i\ssp \TPS)\ess .
\eqno  6.26
$$
Now from 6.23 and 6.25 we get that $A(t)$ and $B(t)$ are $G$-invariant, and  
$G\in \CA$ gives that 
$$
A(t)\ssp \&\ssp B(t)\in \CF[t]\ess .
$$
Since by construction $B(t)$ has distinct roots we can find two polynomials $p(t),q(t)\in \CF[t]$
such that 
$$
p(t)B'(t)+q(t)B(t)\ses 1\ess .
$$
Setting $t=\ggg \TPS$ gives
$$
p(\ggg\TPS)B'(\ggg\TPS)\ses 1
\eqno 6.27
$$
On the other hand we see from 6.26 that if $\ggg=\tau_ih$ with $h\in \TG_\PS$ then 
$$
A(\ggg\TPS)\ses A(\tau_i\TPS)\ses \tau_i\TPP\ssp  B'(\tau_i\TPS)\ess . 
$$
Using 6.23 we may rewrite this as
$$
A(\ggg\TPS)\ses \ggg\TPP\ssp B'(\ggg \TPS)\ess . 
$$
Multiplying both sides by $p(\ggg\TPS)$ and using 6.27 we finally get
$$
\ggg\TPP\ses p(\ggg\TPS)A(\ggg\TPS)
$$
which is 6.24 with
$$
\TT(t)\ses p(t)A(t)\ess .
$$ 
\sa

This Theorem has an immediate consequence which can be helpful in the construction of the Galois group 
of an equation. 
\sa

\noindent{\bol Corollary 6.1} 

{\ita Let $G$ be the Galois group of $\TE_n(t)\in\CF[t]$ and suppose that for
some $\PS\in \CF[\xon]$ we have 
$$
\eqalign{
&(1)\ess\ess S_\PS = G_\PS\ess ,\cr
&(2)\ess\ess \TPS\in \CF \ess .\cr }
$$
Then we must necessarily have }
$$
G\ssp \con\ssp G_\PS\ess .
$$ 
\noindent{\bol Proof}

In view of Theorem 6.2 we need only show that $G_\PS\in \CA(\CF)$. To this end let 
$\PP\in \CF[\xon]$ and let $S_\PP\son G_\PS$. Condition $(1)$ then assures that $\TG_\PP\son \TG_\PS$.
We can thus use Theorem 6.4 and derive that for some $\theta(t)\in\CF[t]$ we have
$\TPP=\theta(\TPS)$. But then condition $(2)$ yields us that $\TPP\in \CF$ as desired. 
\sa

Before we can proceed any further we need to establish the following basic fact
\sas

\noindent{\bol Proposition 6.3}  

{\ita  Let  $B(t)\in\CF[t]$ be a polynomial of degree $k$ which is irreducible in $\CF[t]$, and let 
$$
B(\TPS)\ses 0
$$
for some $\PS\in \CF[\xon]$. Then the values
$$
1\scs \TPS\scs \TPS^2\scs \ldots ,\TPS^{k-1}
\eqno 6.28
$$
form a basis of a vector space $V$ over $\CF$ which is also a field}
\sas

\noindent{\bol Proof}

Suppose that for some $c_o,c_1,c_2,\ldots ,c_{k-1}\in \CF$ not all vanishing we had
$$
c_o+c_1\TPS+c_2\TPS^2+\cdots +c_{k-1}\TPS^{k-1}\ses 0
$$
then the polynomial $R(t)=c_o+c_1t+c_2t^2+\cdots +c_{k-1}t^{k-1}$ would have a root in common
with $B(t)$ and the greatest common divisor of $R(t)$ and $B(t)$ would yield a non trivial 
factorization of $B(t)$ in $\CF[T]$ contradicting the irreducibility of $B(t)$. 
This shows that the elements in 6.28 are independent over $\CF$.

To complete the proof we need to show that every non vanishing element $v$ of $V$ has an inverse in $V$.
Now such an element would be given by a linear combination
$$
v\ses c_o+c_1\TPS+c_2\TPS^2+\cdots +c_{k-1}\TPS^{k-1}\bigsp (\ess c_i\in\CF\ess )
$$
with some $c_i\neq 0$. For the same reasons as above, the polynomial $R(t)=c_o+c_1t+c_2t^2+\cdots +c_{k-1}t^{k-1}$
cannot have any root in common with $B(t)$. Thus we can use the Berlekamp algorithm and construct two polynomials
$ p(t),q(t)\in\CF[t]$ giving
$$
p(t)R(t)+q(t)B(t)\ses 1\ess .
$$
Setting $t=\TPS$ we get  
$$
p(\TPS)\ssp R(\TPS)\ses 1\ess ,
$$
which shows that $ R(\TPS)$ is invertible and that its inverse in $V$ is given by $p(\TPS)$.

\hfill Q.E.D.
\sa

The vector space $V$ will here and after be denoted by $\CF[\TPS]$ and referred to as the {\ita Extension of $\CF$}
by $\TPS$. The integer $k$ giving the dimension of $\CF[\TPS]$ will be called the {\ita degree} of
the extension. We shall also say that  $\CF[\TPS]$ is obtained by  {adjoining}  $\TPS$ to $\CF$. 
The following theorem provides the crucial tools needed in the applications of Galois theory
to the theory of equations. 
\sas

\vbox{
\noindent{\bol Theorem 6.5}

{\ita Let $\PS\in \CF[\xon]$, set
$$
G\ses \tau_1\TG_\PS\sps \tau_2 \TG_\PS + \cdots +\tau_k \TG_\PS
\bigsp (\ess \tau_1=id\ssp )
\eqno 6.29
$$
and let $\TPS_1=\tau_1\TPS,\TPS_2=\tau_2\TPS,\ldots ,\TPS_k=\tau_k\TPS$ denote the conjugates of $\TPS$ in $G$.
Then  
\sas

\item {(i)} The polynomial $B(t)=\prod_{i=1}^k\ssp (t-\TPS_i)$ is irreducible in $\CF[t]$.
\sas

\item {(ii)} By adjoining $\TPS$ to $\CF$ the Galois group of the equation
$\TE_n(t)=0$ is reduced to $\TG_\PS$. 
\sas

\item {(iii)} The Galois group of $B(t)$ is the subgroup  $\Gamma_\TPS$ of $S_k$ corresponding to the action of $G$ on the left
cosets of $\TG_\PS$. In particular we have the isomorphism
$$
\Gamma_\TPS\ssp \cong G\ssp / \ssp \bigcap_{i=1}^k\ssp \tau_i\TG_\PS\tau_i^{-1}\ess .
\eqno 6.30
$$

\item {(iv)} Set $\CF_1=\CF[\TPS]$ and let $Aut_\CF(\CF_1)$ denote the group  of automorphisms of $\CF_1$ 
that leave $\CF$ elementwise fixed. This given, 
$Aut_\CF(\CF_1)$ can be identified with the set $\{\tau_i\TG_\PS  :  \tau_i\TG_\PS=\TG_\PS\tau_i\ssp\}$. In
particular we have the isomorphism 
$$
Aut_\CF(\CF_1)\ssp \cong\Bigl( \sum_{\tau_i\TG_\PS=\TG_\PS\tau_i}\ess \tau_i\TG_\PS\Bigr) /\TG_\PS
\eqno 6.31
$$
}
}
\noindent{\bol Proof}

\item {\ita Proof of (i)}

Suppose that a polynomial
$B_1(t)\in \CF[t]$ divides $B(t)$ and shares the root $\tau_i\TPS$ with $B(t)$. Then the polynomial
$B_1(\tau_i\PS)\in \CF[\xon]$ has the value $B_1(\tau_i\TPS)=0\in \CF$ and $G\in \CB$ gives that
$S_{B_1(\tau_i\PS)}\son G$. In other words we must have $B_1(\ggg\tau_i\TPS)=0$ for all $\ggg\in G$.
Since the action of $G$ on $\TPS_1,\TPS_2,\ldots ,\TPS_k$ is transitive, we see that
$B_1(t)$ has to vanish at all the roots of $B(t)$,  which forces $B_1(t)=B(t)$. Thus
$B(t)$ can't have a proper factor in $\CF[t]$. 
\sas

\item {\ita Proof of (ii)}

Since $B(\TPS)=0$ we can use Proposition 6.3 to construct the field $\CF_1=\CF[\TPS]$.
This given  we need only verify that
$$
\CA(\CF[\TPS])\cap\CB(\CF[\TPS])\ses \{\TG_\PS\}\ess .
$$
Now this is immediate. In fact, if for some $\PP\in \CF[\xon]$ we have $S_\PP\son \TG_\PS$
then we must have $\TG_\PP\son\TG_\PS$ as well and Theorem 6.4 gives that for some $\TT(t)\in \CF[t]$
$$
\TPP=\TT(\TPS)\ssp \in \CF[\PS]\ess .
$$
Thus $\TG_\PS$ is in $\CA(\CF[\TPS])$. 

Conversely let  $\PP\in \CF[\xon]$ and $\TPP\in \CF[\TPS]$. This means that for some $\TT(t)\in \CF[t]$ we have
$$
\TPP=\TT(\TPS)\ess .
$$
In other words for the polynomial $\Xi(x)= \PP(x)-\TT(\PS(x))\in \CF[\xon]$ we have
$$
{\tilde \Xi}= 0\ssp \in \CF
$$
so from $G\in \CB(\CF)$ we get that we must have 
$$
\ggg\TPP-\TT(\ggg\TPS)=\ggg {\tilde \Xi}\ses 0\bigsp (\ess \forall\ess\ess \ggg\in G\ess )
$$
However, for $\ggg\in \TG_\PS$ this yields that
$$
\ggg\TPP\ses \TT(\TPS)\ses \TPP\ess .
$$
In other words $\TPP\in \CF[\PS]$ implies that  $S_\PP\son \TG_\PS$, which gives $\TG_\PS\in \CB(\CF[\TPS])$
as desired. This completes the proof of $(ii)$.
\sas

\item {\ita Proof of (iii)} 

\def \GTPS {\Gamma_\TPS} 

Let $\Gamma_\TPS$ be the image of $G$ in the symmetric group $S_k$ given by the permutation action of $G$
on the left cosets of $\TG_\PS$. We want to show that
$$
 \CA_{B}(\CF)\cap\CB_B(\CF)\ses \{\Gamma_\TPS\}\ess .
\eqno 6.31
$$
To this end let $\PP\in \CF[y_1,y_2,\ldots ,y_k]$ and suppose that for any $\ggg\in \GTPS$ we have
$$
\PP(\TPS_{\ggg_1},\TPS_{\ggg_2}, \cdots ,\TPS_{\ggg_k})\ses \PP(\TPS_1,\TPS_2,\ldots ,\TPS_k)\ess .
\eqno 6.32
$$  
Since each $g\in G$ induces a permutation of $\TPS_1,\TPS_2,\ldots ,\TPS_k$ by an element
$\ggg\in \GTPS$, 6.32 implies that
$$
S_{\PP(\PS_1,\PS_2,\ldots ,\PS_k)}\ssp \son \ssp  G\ess .
$$
But then $G\in \CA_{\TE_n}(\CF)$ gives that $\PP(\TPS_1,\TPS_2,\ldots ,\TPS_k)\in \CF$. Thus $\GTPS\in \CA_{B}(\CF)$.

Conversely, say $\PP(\TPS_1,\TPS_2,\ldots ,\TPS_k)=a \in \CF$. Then $G\in \CB_{\CE_n}(\CF)$ gives
that $g\PP(\TPS_1,\TPS_2,\ldots ,\TPS_k)=\PP(\TPS_1,\TPS_2,\ldots ,\TPS_k)$ for all $g\in G$. But if
$\ggg=\ggg(g)$ is the image of $g$ in $S_k$, this simply says that  
$$
\PP(\TPS_{\ggg_1},\TPS_{\ggg_2}, \cdots ,\TPS_{\ggg_1})\ses \PP(\TPS_1,\TPS_2,\ldots ,\TPS_k)\ess .
$$
This gives $\GTPS\in \CB_B(\CF)$. To complete the proof of $(iii)$ we need only observe that the 
Galois stabilizer of a conjugate $\TPS_i$ is simply the conjugate subgroup $\tau_i\TG_\PS\tau_i^{-1}$, 
thus the only elements of $g$ that leave invariant all the conjugates of $\TPS$ are those that belong to 
the intersection  
$$
\cap_{i=1}^k\ssp \tau_i\TG_\PS\tau_i^{-1}\ess .
$$
This yields 6.30.
\sas

\item {\ita Proof of (iv)} 

We assume here that $\TE_n(t)\in \CF[t]$ as before and that $\CE=\CF[\aon]$. 
Note that since every element of $\CF_1$ is of the form 
$$
\theta(\TPS)\ses c_o+c_1\TPS+\cdots +c_{k-1}\TPS^{k-1}\bigsp (\ess c_i\in \CF\ess )
\eqno 6.33
$$
with $k$ the index of $\TG_\PS$ in $G$, to find the image  by a  $g\in Aut_\CF(\CF_1)$
of any element of $\CF_1$  we need only know $g\TPS$. This is because, every $c_i$ remaining
unchanged by $g$ we must necessarily have $g\theta(\TPS)=\theta(g\TPS)$. Moreover, since
$B(t)=\prod_{i=1}^k(t-\TPS_i)\in\CF[t]$, the identity $B(\TPS)=0$ forces $B(g\TPS)=0$
as well. In particular, we deduce that $g\TPS=\TPS_i=\tau_i\TPS$ for some
$i=1,2,\ldots ,k$. In addition $g\TPS\in\CF_1$ yields that we must have $\TPS_i=\tau_i\TPS=\theta(\TPS)$
with $\theta(t)\in \CF$. But this gives that $\TG_\PS\con \TG_{\PS_i}$ and since 
$\TG_{\PS_i}=\tau_i\TG_\PS\tau_i^{-1}$ the latter inclusion can hold true if and only if
$$
\TG_{\PS_i}=\TG_\PS\ess .
$$
Thus we see that the elements of $Aut_\CF(\CF_1)$ can be simply identified with the left cosets
$\tau_i\TG_\PS$ such that $\tau_i\TG_\PS=\TG_\PS\tau_i$. This gives 6.31 and completes the proof of the Theorem.
\sa

We should note that the argument used in the proof of part $(i)$ of this Theorem leads to 
the following basic property of the Galois group of an equation.  
\sa

\noindent{\bol Proposition 6.4}

{\ita The polynomial $B(t)=\prod_{i=1}^n(t-\aaa_i)\in\CF[t]$ is irreducible in $\CF[t]$ if and only 
if its Galois group $G=G_B(\CF)$ acts transitively on $\aon$. }
\sas

\noindent{\bol Proof}

Let $B_1(t)\in \CF[t]$ be an irreducible factor of $B(t)$. If $\aaa_i$ is any of the roots of $B_1(t)$ then
$B_1(\aaa_i)=0\in \CF$ and the fact that $G\in\CA_B(\CF)$ gives that $B_1(\aaa_{\gamma_i})=0$
for all $\gamma\in G$. But then the transitivity of $G$ yields that $B_1(t)$ can't be a proper factor of $B(t)$.
Conversely suppose that $G$ is intransitive. This means that the orbit of any of the roots, say $\aaa_1$
can't consist of all the roots of $B(t)$. Denote this orbit by $Orb(\aaa_1)$ and set
$$
B_1(t)=\prod_{\aaa_i\in Orb(\aaa_1)}(t-\aaa_i)\ess .
$$
Since $B_1(t)$, by construction is invariant under $G$ then $G\in\CA_B(\CF)$ gives that
$B_1(t)\in \CF[t]$. Moreover, also by construction, $B_1(t)$ does not contain all then roots of $B(t)$.
Thus $B_1(t)$ is necessarily a proper factor of $B(t)$ and the latter must therefore be reducible
in $\CF[t]$.   
\sa

\noindent{\bol Remark 6.3}

We should note that parts $(i)$, $(ii)$ and $(iii)$ of Theorem 6.5 are Galois' fundamental breakthroughs
in the theory of equations. In the original Galois context groups came first and fields were only accessories.
In later interpretations and additions to Galois' work started by Kronecker [], brought to completion by
Dedekind [] and Weber and popularized by E. Artin [],[], this viewpoint has been reversed and
Galois theory was made to become part of the theory of fields. In particular,   
part $(iv)$ of Theorem 6.5 in its present interpretation is a later addition.

Since our presentation of Galois theory differs from most available textbooks on the subject, perhaps a few
words might be needed to connect this writing to present day literature. For instance in Artin's monograph []
an extension field $\CF_1$ of a field $\CF$ is called a {\ita normal} extension if {\ita the group of 
automorphisms of $\CF_1$ that leave $\CF$  fixed has $\CF$ for its fixed field}. We should note that in []
a field $\CF$ is said to be {\ita fixed } by an automorphism $g$ if $g$ fixes {\ita every} element of $\CF$.
We see then that one of our fields $\CF_1=\CF[\PS]$ is a normal extension of $\CF$ if and only if the only
elements of $\CF_1$ that remain fixed under every element of $Aut_{\CF}(\CF_1)$ (as given by 6.31)
are the elements of $\CF$ itself. 

This given, we can easily convert Theorem 6.5 into the collection of results that in [] is referred to
as the {\ita Fundamental Theorem of Galois Theory}.
\sas

\noindent{\bol Theorem 6.6}

{\ita Let $G$ be the Galois group of $\TE_n(t)=\prod_{i=1}^n(t-\aa_i)$ in $\CF$ and
let $\CE=\CF[\aon]$. Then
\sas

\item{(i)} Each subgroup $G_1\con G$ is the Galois group of $E_n(t)$ with respect
to an intermediate field $\CF_1=\CF[\TPS]$. Different groups $G_1,G_2$  corresponding to
different fields $\CF_1,\CF_2$. 
\sas

\item{(ii)} The subgroup $G_1$ is a normal subgroup of $G$ if and only if the corresponding field
$\CF_1$ is a normal extension of $\CF$. In that case the group of automorphisms of $\CF_1$ that
leave $\CF$ fixed  is isomorphic to the quotient group $G/G_1$. 
\sa
 
\item{(iii)} For each $G_1\con G$ the dimension of $\CF_1$ over $\CF$ is $|G|/|G_1|$ and the 
dimension of $\CE$ over $\CF_1$ is $|G_1|$.}
\sa

\noindent{\bol Proof}

\item{\ita proof of (i)}

We have seen (Theorem 6.1) that given any subgroup $G_1\con G$  we can find $\PS\in \CF[\xon]$ such that $\TG_\PS=G_1$.
From Theorem 6.5 we get that the Galois group of $\TE_n$ in  $\CF_1=\CF[\PS]$ is $G_1$. If for $G_1,G_2\con G$
we have $\TG_{\PS_1}=G_1$, $\TG_{\PS_2}=G_2$ then  $\CF[\TPS_1]\con\CF[\TPS_2]$ gives $\TPS_1\theta(\TPS_2)$
with $\theta(t)\in \CF[t]$ thus also $\TG_{\PS_2}\con \TG_{\PS_1}$ so the equality $\CF[\PS_1]=\CF[\PS_2]$ 
forces the equality $G_1=\TG_{\PS_1}=\TG_{\PS_2}=G_2$. This proves $(i)$. 
\sas

\item{\ita proof of (ii)}

Suppose that $G_1=\TG_\PS$ is not a normal subgroup of $G$. Then by relabeling the elements $\tau_i$ 
appearing in 6.29 so that 6.34 may be rewritten as
$$
Aut_\CF(\CF_1)\ses \sum_{i=1}^s \ssp \tau_i\TG_\PS\ess \bigsp (\ess for \ess\ess s<k\ess )
$$
Now note that the polynomial  
$$
B_1(t)\ses \prod_{i=1}^s (t-\TPS_i)\in \CF_1[t]
$$
cannot be in $\CF[t]$. This is clear since otherwise the irreducibility  of $B(t)= \prod_{i=1}^k (t-\TPS_i)$
(Theorem 6.5 $(i)$) would be contraddicted. In particular one of the elementary symmetric functions
$e_i(y_1,y_2,\ldots ,y_s)$ must take a value $e_i(\TPS_1,\TPS_2,\ldots ,\TPS_s)$ not in $\CF$. However,
since the latter is invariant under permutations of $\TPS_1,\TPS_2,\ldots ,\TPS_s$ we see (from 6.31) that we
have a element of $\CF_1$ that remains invariant under all elements of $Aut_\CF(\CF_1)$. So
$\CF_1$ in this case is not a normal extension of $\CF$. 
Conversely suppose that $G_1=\TG_\PS$ is a normal subgroup of $G$ then 6.31 gives that
$$
Aut_\CF(\CF_1)\ses \Bigl( \sum_{i=1}^k\ess \tau_i\TG_\PS\Bigr) /\TG_\PS \ses G/G_1\ess .
\eqno 6.34
$$
Since every  element of $\CF_1$ is already fixed by $\TG_\PS$, we see from 6.34 that an element $f\in\CF_1$ is
fixed by $Aut_\CF(\CF_1)$ if and only if it is fixed by every element of $G$, but then $G\in\CA_{\TE_n}(\CF)$ gives that 
$f$ must lie in $\CF$. This gives that $\CF$ is the fixed field of $Aut_\CF(\CF_1)$ and completes the proof of $(ii)$.
\sas

\item{\ita proof of (iii)}

Given that $G_1=\TG_\PS$ and given that we have 6.29, Proposition 6.3 combined with $(i)$ of Theorem 6.5 
yields that the dimension of $\CF_1$ as a vector space over $\CF$ is precisely given by $k=|G|/|G_1$.
It develops that the last assertion of the Theorem is an immediate consequence of part $(ii)$ 
of Theorem 6.5. In fact, to construct the Galois group of $\TE_n$ relative to $\CF_1$ we can also
resort to the original definition. That is we break up the polynomial $\TP(t)$ given in 6.11 into
its irreducible factors in $\CF_1[t]$ and take the collection of permutations which give the irreducible factor $P_{1,1}(t)$
that has $\tv$ as a root.   
\footnote {(*)}{Note that $P_{1,1}(t)$ must also be an irreducible factor of the polynomial $\TP_1(t)$
which gave us the Galois group $G$ of $\TE_n$ in $\CF$.}

However, part $(ii)$ of Theorem 6.5 gives that this procedure must deliver
the polynomial
$$
P_{1,1}(t)\ses \prod_{\tau\in \TG_\PS}\ssp (t-\tau \tv)\ess .
\eqno 6.35
$$
Since $P_{1,1}(t)$ is in $\CF_1[t]$ we can use Proposition 6.3 and deduce that the extension
$\CF_1[\tv]$ must be of dimension $|\TG_\PS|$ as a vector space over $\CF_1$. But then $(iii)$
follows from the fact that $\CF_1[\tv]$ and $\CF[\aon]$ are one and the same.  
In fact, we trivially have $\CF_1[\tv]\con\CF[\aon]$ and the reverse containement $\CF_1[\tv]\son\CF[\aon]$ 
follows from Theorem 6.1. This completes our proof. 
\sa

\noindent{\bol Remark 6.4}

We should note that under the definition of some texts (see [] ,[]) one of our extensions $\CF_1=\CF[\TPS]$ 
would be called $normal$ if and only if any irreducible polynomial $Q(t)\in\CF[t]$ that has a root in $\CF_1$
has all the other roots in $\CF_1$. Now we can easily show that, again this happens if and only if $\TG_\PS$ 
is a normal subgroup of the Galois group $G=G_{\TE_n(t)}(\CF)$. Let us use the same notation as in the
proof of Theorem 6.5 and let us set $Q(t)=\prod_{i=1}^m(t-\bbb_i)$. Note first that if the ``normality'' condition 
in [] and [] is satisfied then one of the polynomials all whose roots,
by this condition, would have to lie in $\CF_1$ is the polynomial $B(t)$ itself.  But as we have seen this
is equivalent to $\TG_\PS$ being a normal subgroup of $G$. 

Conversely suppose that  $\TG_\PS$ is a normal subgroup of $G=G_{\TE_n(t)}(\CF)$. Then  $\bbb_1\in\CF_1$, 
simply means that $\bbb_1=\theta(\TPS)$ with $\theta(t)\in \CF[t]$. In other words we have
$$
\bbb_1=\PP(\aon)=\theta(\PS(\aon))\ess .
$$
But then we must also have
$$
\ggg\bbb_1=\ggg\TPP=\theta(\ggg\TPS).
$$
Letting $\Delta$ denote a set of representatives for the left
cosets of $\TG_\PP$ in $G$ construct the polynomial
$$
Q_1(t)\ses \prod_{\tau\in \Delta}(t-\tau\TPP)\in \CF[t]
$$ 
Since $Q_1(t)$t shares a root with $Q(t)$, the irreducibility of the latter forces 
all the roots of $Q(t)$ to be roots of $Q_1(t)$. In other words, every root $\bbb_i$ has
an expression of the form
$$ 
\bbb_i= \tau_{j_i}\TPP= \theta(\tau_{j_i}\TPS )\ess . 
$$
Now, since the normality of $\TG_{\PS}$ forces all the conjugates of $\TPS$ in $G$ to belong
to $\CF[\TPS]$ we must conclude that each $\bbb_i\in\CF[\TPS]$ as desired. 

We should note that $Q(t)$ being a factor of $Q_1(t)$ yields that $m\leq k=degree\ssp B(t)$.
\sap
\def \om {\omega}
\def\CF {{\cal F}}
\def\CE {{\cal E}}
\def\CA {{\cal A}}
\def\CB {{\cal B}}

\def\RFX {Rat[{\cal F};x_1,\ldots ,x_n]}
\def\PFX {{\cal F}[x_1,\ldots ,x_n]}
\def \xsig {x_{\sig_1},x_{\sig_2},\ldots ,x_{\sig_n}}
\def \PP {\Phi}
\def \PS {\Psi}
\def \SYM {Sym[\CF;\xon]}
\def \RSYM {Ratsym[\CF;\xon]}
\def \part {\vdash}
\def \DD {\Delta}
\def \omi {\hskip -.01truein\ ^{(i)}}
\def \scos {\ssp :\ssp}
\def \bup {\hskip -.01truein}
\def \om {\omega}
\def \eon {e_1,e_2,\ldots ,e_n}
\def \SQR {{\root 2 \of R}}
\def \sh {{\scriptstyle{1\over 2}}}
\def \shf {{\scriptstyle{1\over 4}}}
\def \AO {{\root 3 \of {-\sh q+ \SQR}}}
\def \BO {{\root 3 \of {-\sh q- \SQR}}}
\def \aon {\aaa_1, \aaa_2,\ldots ,\aaa_n}

\def \TPP {{\tilde \PP}} 
\def \TE {{\tilde E}} 
\def \TT {\Theta} 
\def \TP {{\tilde P}} 
\def \TQ {{\tilde Q}}
\def \tv  {{\tilde v}}
\def \TPS  {{\tilde \PS}}
\def \TG {{\tilde G}}
\def \RA {\ssp \Longrightarrow\ssp}

\def \nor {\triangleleft}
\def \son {\supseteq}
\def \res {\TP_1(t)}
\def \fora {\bigsp (\ess \forall \ess\ess \ggg\in G\ess )}
\def \ttt {\theta}
\def \CFT {\CF[t]}
\def \EN {{\TE_n(t)}}
\noindent{\bol 7. Solving cyclic equations}
\sa

In these notes we shall say that an $n^{th}$ degree  polynomial $\TE_n(t)=\prod_{i=1}^n(t-\aaa_i)\in \CFT$ and the corresponding
equation $\TE_n(t)=0$ is $cyclic$ in $\CF[t]$ if its Galois group $G=G_\EN (\CF)$ is the cyclic group on $n$ letters.
More precisely, $\TE_n(t)$ is cyclic with respect to $\CF[t]$ if and only if by a suitable labeling of the roots $\aon$
we have
$$
G=G(\CF)=\{id,\ggg,\ggg^2,\ldots ,\ggg^{n-1}\}
\eqno 7.1
$$
with
$$
\ggg\aaa_i\ses \aaa_{i+1}
\bigsp  for \ess\ssp i=1,2,\ldots ,n\ess\ess\ess (\aaa_{n+1}=\aaa_1)
\eqno 7.2
$$
\sas
We should note that in our definition of cyclicity we implicitely assume that $E_n(t)$ has
distinct roots. This given, the following basic fact is helpful in establishing cyclicity.
\sas

\noindent{\bol Theorem 7.1}

{\ita $\TE_n(t)$ is cyclic with respect to $\CF[t]$ if and only if
\sas

\item {(1)} It is irreducible in $\CF[t]$,
\sas

\item {(2)} We have a polynomial $\ttt(t)\in \CF[t]$ and a labeling $\aon$ of the 
roots of $\TE_n(t)$ such that}
$$
\aaa_{i+1}\ses \ttt(\aaa_i)
\bigsp  for \ess\ssp i=1,2,\ldots ,n\ess\ess\ess (\aaa_{n+1}=\aaa_1)
\eqno 7.3
$$
\noindent{\bol Proof}

Suppose $\EN$ is cyclic. Then under the labeling which gives 7.1 and 7.2 set
$$
Q(t)\ses \sum_{s=1}^n \ssp \aaa_{s+1}\ssp \prod_{j=1}^{(s)}(t-\aaa_j)\ess .
\bigsp  (\ssp \aaa_{n+1}=\aaa_1\ssp)
$$
Note that for $i=1,2,\ldots ,n$ we have
$$
Q(\aaa_i)\ses \aaa_{i+1}\ssp \TE_n'(\aaa_i)
\eqno 7.4
$$
Since the roots $\aon$ are supposed to be distinct $\TE_n(t)$ and  its derivative $\TE_n'(t)$ have no roots
in common. Thus using the Berlekamp algorithm we can construct two polynomials $A(t),B(t)\in \CFT$ such that
$A(t)\TE_n(t)+B(t)\TE_n'(t)=1$. Since setting $t=\aaa_i$ yields $B(\aaa_i)\TE_n'(\aaa_i)=1$, we see that we can rewrite
7.4 in the form
$$
\ttt(\aaa_i) \ses \aaa_{i+1} 
$$
with
$$
\ttt(t)\ses B(t)Q(t)\in \CFT\ess .
$$
This proves property $(2)$. To show $(1)$ we need only observe that, in view of the transitivity of the action in 7.2,
the irreducibility of $\EN$ follows from Proposition 6.4.
\sas

Suppose now that $\EN$ satisfies $(1)$ and $(2)$. Then, using the labeling that gives 7.3, we may define $\ggg$
as the circular permutation that gives $\ggg\ssp \aaa_i=\aaa_{i+1}$. Let us also recursively define the
polynomials $\ttt_{i}(t)\in\CFT$ by setting $\ttt_{o}=t$ and
$$
\ttt_{i+1}(t)\ses \ttt(\ttt_{i}(t))\bigsp (\ess i=1,2,\ldots ,n-1\ess )\ess .
$$  
Let as also set for any $\PP(\xon)\in \CF[\xon]$ 
$$
R_\PP(t)\ses \PP\bigl(\ttt_o(t),\ttt_{1}(t),\ttt_{2}(t),\ldots,\ttt_{n-1}(t)\bigr)\ess .
\eqno 7.5
$$
Now let $G=\{id,\ggg,\ggg^2,\cdots ,\ggg^{n-1}\}$ and suppose that 
$$
S_\PP\son G
\eqno 7.6
$$
Using 7.2 and 7.5 we may rewrite this as
$$
R_\PP(\aaa_1)=R_\PP(\aaa_2)=\cdots =R_\PP(\aaa_n)\ess .
\eqno 7.7
$$
In particular, we must have
$$
\PP(\aon)=R_\PP(\aaa_1)\ses {1\over n}  \bigl(R_\PP(\aaa_1)+R_\PP(\aaa_2)+\cdots +R_\PP(\aaa_n)\bigr)
$$
However since 
$$
{1\over n}  \bigl(R_\PP(x_1)+R_\PP(x_2)+\cdots +R_\PP(x_n)\bigr)\in Sym[\CF;\xon]
$$
its value at $\aon$ must be expressible as a polynomial in the coefficients of $\EN$, This gives
that
$$
\PP(\aon)=a\in \CF
\eqno 7.8
$$
and establishes that $G\in\CA_{\EN}(\CF)$. To finish the proof we need only show that
we also have 
$$
G\in\CB_{\EN}(\CF)\ess .
\eqno 7.9
$$
To this end suppose that 7.8 hold true for some $\PP\in \CF[\xon]$.
Using 7.5 we may translate this property into the statement that the polynomial
$$
R_\PP(t)-a\in \CFT
$$
vanishes for $t=\aaa_1$. Now under (1) $\EN$ is irreducible and this forces it to be a factor of
$R_\PP(t)-a$. In other words we must have
$$
R_\PP(\aaa_i)=a\bigsp (\ssp for\ess\ess i=1,2,\ldots ,n\ess )\ess .
$$
This shows that $S_\PP\son G$ and that 7.9 holds true as desired.
\sa

Note that although cyclicity is a field dependent property, here and in the following we shall drop the
appendage ``$in\ess\CF[t]$''  or ``{\ita with respect to $\CF[t]$}'' in all cases in which the identity of
the base field $\CF$ is clear from the context.
\sa

Theorem 7.1 has two immediate important applications:
\sa

\noindent{\bol Corollary 7.1}

{\ita If a binomial equation $t^p-a=0$ with $0\neq a\in \CF$ and $p$ a prime is irreducible in $\CF[t]$, then it 
is cyclic if and only if $\CF$ contains a primitive $p^{th}$-root of unity.}
\sas

\noindent{\bol Proof}

Let $w\in \CF$ be a primitive $p^{th}$-root of unity and let $\aaa$ be a root of $t^p-a$. Then since for
no $s<p$ we may have $w^s=1$ 
\footnote {(*)}{the smallest such $s$ would have to divide $p$}
the powers $w^i$ (for $i=1,..,p-1$) are all distinct. This gives that the roots of $\ess t^p-a\ess $ are simply
$$
\aaa_1=\aaa\scs \aaa_2=w\aaa\scs\aaa_3=w^2\aaa\scs \ldots ,\scs \aaa_n=w^{n-1}\aaa
$$
and so we have 7.3 with $\ttt(t)=wt$. Thus the cyclicity of $t^p-a$ follows from
Theorem 7.1. 

Conversely suppose that $t^p-a$ is cyclic. Let $\aaa_1,\aaa_2,\ldots ,\aaa_p$ be the labeling of
its roots that gives 7.2 so that its Galois group is $G=\{id,\ggg,\ldots ,\ggg^{p-1}\}$. 
Set $w=\aaa_2/\aaa_1=\aaa_2\aaa_1^{p-1}/a$. Clearly we must have $w^p=1$ and we can't have $w=1$ since that
would contradict the irreducibility of $t^p-a$. This given, since the 
the elements $w^{i-1}\aaa_1$ are distinct and all satisfy $t^p-a=0$ they must be a permutation 
of $\aaa_1,\aaa_2,\ldots ,\aaa_p$. Thus
$$
t^p-a\ses \prod_{i=1}^p\ssp(t-w^{i-1}\aaa_1)\ess .
$$
We thus deduce that 
$$
\aaa_i=w^{h_i}\aaa\bigsp (\ess i=1,\ldots p\ess )
$$
where  $h_1,h_2,\ldots, h_p$ is a permutation of $0,1,\ldots ,p-1$. From this we derive that $\ggg
w=\aaa_3/\aaa_2=w^{h_3-1}$. Since $\aaa_3\neq \aaa_2 $ we see that $h_3\neq 1 \ess mod\ess p$. 
Setting $h_3-1=e$ we get $\ggg^kw=w^{e^k}$. But then a theorem of Euler gives that $s^{p-1}\cong 1$ mod $p$.
In particular we get that $\ggg^{p-1}w=w$. Now, $\ggg^{p-1}=\ggg^{-1}$ and $\ggg^{-1}$ generates $G$ as well as
$\ggg$. Thus $w$ is invariant under $G$ and must necessarily belong to $\CF$ as we wanted to show.
\sa
\def \CQ {{\cal Q}}

\vbox{
\noindent{\bol Corollary 7.2}

{\ita The cyclotomic polynomial 
$$
\PP_p(t)= 1+t+t^2+\cdots +t^{p-1}\bigsp (for\ess p\ess a\ess prime\ess)
\eqno 7.10
$$
is cyclic with respect to the field $\CQ$ of rational numbers}
\sas

\noindent{\bol Proof}
}

We start by showing that $\PP_p$ is irreducible in $\CQ [t]$. Suppose not. Then by Gauss theorem
we will have a non trivial factorization 
$$
\PP_p(t)\ses A(t)B(t)
\eqno 7.11
$$
where $A(t)$ and $B(t)$ are  both monic polynomials with integer coefficients. Setting $t=1$
gives
$$
p=\PP_p(1)\ses A(1)B(1)\ess .
$$
So one of $A(1),B(1)$ must be $\pm 1$. Say it is  $A(1)$. Now let $\aaa$ be a root of $\PP_p$.
Since $\aaa^p=1$ we can't have $\aaa^s=1$ for any $1\leq s\leq p-1$. In particular,
for such an $s$, the elements 
$$
\aaa^s\scs\aaa^{ 2s}\scs\ldots \scs\aaa^{(p-1)s}
$$
are just a rearrangement of the roots of $\PP_p$. Since at least one of them is a root of $A(t)$, we are forced to
conclude that
$$
A(\aaa^s)A(\aaa^{ 2s})\cdots A(\aaa^{ (p-1)s})\ses 0\ess .\bigsp (\ess\forall\ess 1\leq s\leq p-1\ess)
$$
This yields that the polynomial
$$
R(t)=A(t)A(t^{2})\cdots A(t^{(p-1)})
$$ 
must be divisible by $\PP_p$. Since $R(t)$ is monic with integer coefficients, again by Gauss
theorem, we shall have the factorization
$$
R(t)=(1+t+t^2+\cdots +t^{p-1})R_1(t)
$$
with $R_1(t)$ also monic with integer coefficients. Now setting $t=1$ we are forced to the impossible
conclusion that
$$
\pm 1 \ses p\ssp R_1(1)\ess .
$$
Thus $\PP_p$ must be irreducible as asserted. 
\sas

We know from number theory that for any prime $p$ we can find a primitive exponent $e\in [1,p-1]$ 
which has the propery that the integers 
$$
e\scs e^2\scs \ldots e^{p-2}
$$
are distinct and, in fact, give (modulo $p$) a permutation of
$$
1\scs 2 \scs \ldots , p-1
$$
Choosing one such exponent, we derive that if $\aaa$ is any of the roots of $\PP_p$ then the powers
$$
\aaa\scs \aaa^e\scs \ldots \aaa^{e^{p-2}}
$$
give back again all the roots of $\PP_p$. But this means that if we label the roots of $\PP_p$
by setting $\aaa_i=\aaa^{e^{i-1}}$ we shall have 7.3 with
$$
\theta(t)\ses t^e\in \CQ [t]\ess .
$$
This completes our proof.
\sas

There is a further property of the binomial equation with prime exponent we need to  know here.
\sa

\noindent {\bol Proposition 7.1}

{\ita The binomial $t^p-a$ ($a\in \CF$) is reducible in $\CF[t]$ if and only if $a$ is the
$p^{th}$ power of an element of $\CF$.}
\sas

\noindent{\bol proof}

Suppose that $t^p-a$ is reducible in $\CF[t]$ and let 
$$
t^p-a\ses A(t)B(t)\bigsp (\ess A(t),B(t)\in \CF[t]\ess )
$$
be a non trivial factorization. We may then write
$$
A(t)=\prod_{i=1}^s(t-w^{h_i}\aaa)
$$
where  $s<p$, $\aaa$ is any root of $t^p-a$, and $w$ can be chosen to be the ratio of any two roots.
However, $A(t)\in \CF[t]$ implies in particular that its constant term $c$ is in $\CF$.
In other words 
$$
\aaa^s w^{({h_1}+{h_2}+\cdots +{h_s})}=c \ess \in \CF
\eqno 7.12
$$
Since $1\leq s<p$ we can find integers $h,k$ such that $hs=1+kp$. Raising both sides of 7.12 to the power $h$
yields 
$$
\aaa\ssp  a^k w^{h'}=c^h\ess ,
$$ 
for a suitable integer $0\leq h'\leq p-1$. Thus one of the roots of $t^p-a$, namely $b=\aaa  w^{h'}$
lies in $\CF$ and $a$ must necessarily the $p^{th}$ power of an element of $\CF$ as we asserted.

The converse is entirely trivial since when $a=b^p$ with $b\in \CF$ we have the factorization
$$
t^p-a=(t-b)(b^{p-1}+b^{p-2}t+\cdots +bt^{p-2}+t^{p-1})\ess .
$$
\sa

\def\CP {{\cal P}}

To proceed any further in this section we need to update the meaning of
{\ita solving by radicals} in the Galois setting. It is natural to
assume that in this setting {\ita root extraction} should simply mean 
extending a given field $\CP$, by the {\ita adjunction} of a root of
a binomial equation
$$
t^n-a=0\ess .\bigsp (\ess a\in \CP\ess)
$$ 
This given, solving by radicals the equation $\TE_n(t)=0$, in the
Galois setting means constructing the roots  of $\TE_n(t)\in \CF[t]$
by a sequence of extensions
$$
\CF_{k-1}\ssp \rightarrow\ssp \CF_{k}\ses \CF_{k-1}[\xi_{k}]\ess\ess\ess \bigsp (\ess k=0,1,2,..,k_o\ess )
\eqno 7.13
$$
with
$$
\eqalign{
& (1)\ess\ess \xi_{k}^{p_k}\in \CF_{k-1}\cr  
&(2)\ess\ess \CF_o=\CF\ess ,\cr
&(3) \ess \aon\in \CF_{k_o}\ess .\cr }
\eqno 7.14
$$
In particular, this will enable us to construct a formula,
say for $\aaa_1$, which will be of the form
$$
\aaa_1 \ses 
\cdots \sps {\root p \of 
{\cdots 
    {\root q\of 
         {\cdots + 
            {\root r \of \cdots }+ {\root s \of \cdots }
          \cdots }
    }
}
}
\eqno 7.15
$$
Where the radicals appearing in it will be appropriately chosen solutions
of the equations 
$$
t^p_k - a_k\ses 0 \ess\ess\ess\ess {\rm with}\ess\ess a_k\in \CF_{k-1}
\eqno 7.16
$$
If we chose to  represent the element $\xi_{k}$ appearing in 7.14 $(1)$
by the symbol $\root p_k \of a_k$, we are confronted with the ambiguity
resulting from the multivalued nature of the symbol ``$\root p \of ...$''.
For instance if, at the $k^{th}$ step of the extension process, $\xi_k$ 
is to be a primitive $6^{th}$-root of unity. It would be better to represent
$\xi_k$ by  $-{1\over 2}+{{\root 2 \of -3}\over 2}$ rather by $\root 6 \of 1$.
This is because the expression $-{1\over 2}+{{\root 2 \of -3}\over 2}$
yields only these two primitive roots as we specialize $\root 2 \of -3$ 
to the two conjugate roots of the equation $x^2+3=0$. Clearly, it is in the nature
of the problem that any formula we may construct for $\aaa_1$ should 
have multivalued symbols appearing in it. For $\aaa_1$ itself is in essence an $n$-valued
function. Indeed, labeling the roots of $\TE_n$  $\aon$ is an artificial device, and
``$\aaa_1$'' should really represent only a {\ita generic }  root of $\TE_n$. 
Now, we have seen already in the Lagrange setting, that there is no loss in requiring
that in the successive radicals $\root p_k \of a_k$ the exponent $p_k$ should be a prime number. 
Now if the equation $t^{p_k}-a_k=0$ is {\ita cyclic} in $\CF_{k-1}$ the adjonction of any
one of its roots to $\CF_{k-1}$ will result in the same field $\CF_k$. So writing
$\CF_k=\CF_{k-1}[{\root p_k \of a_k}]$ does not produce any ambiguity as far $\CF_k$ is concerned.
Moreover, we see that if we require that each symbol $\xi_k={\root p_k \of a_k}$ appearing in 7.14
represents the same root of $t^{p_k}-a_k=0$ and if we let each of these symbols in turn and
independently describe each of the other solutions of the corresponding equation $t^{p_k}-a_k=0$
then $\aaa_1$ as expressed by 7.14 will represent a $p_1p_2\cdots p_{k_o}$-valued
function. So that if $p_1p_2\cdots p_{k_o}=n$ then  7.14 by this process
will deliver each of the roots of $\TE_n$. Here and after a formula 7.15
satisfying these requirements will be called a {\ita tight } formula 
and the symbols ${\root p_k \of a_k})$ appearing in it will likewise be called  {\ita tight radicals}.
More generally we shall refer to an adjunction $\CP\rightarrow \CP[\xi={\root p \of a}]$ as a 
{\ita tight radical extraction} if and only if $t^p-a$ is cyclic in $\CP$. All other formulas
and root extractions will be referred to as {\ita loose}. It is not difficult to verify that
the formulas for the roots of the cubic and quartic which can be obtained by the process
given in section 5 are, in fact, tight. However as we pointed out, the process in section 5 has the additional
property that at each step the element adjoined (roots of unity apart), is a root
of a polynomial function of the roots $\xon$. This brings us to define 7.15 as a {\ita natural}
formula if every $\xi_k$ is also a polynomial in $\aon$ with coefficients in $\CF$.  
In this terminology, we can say that Lagrange showed that (as long as $\CF$ contains all
the needed roots of unity) the roots of the general cubic and quartic 
may be given by a tight natural formula, Ruffini showed that the roots of the general
quintic have no loose natural formulas and Abel showed that they may not even be 
solved by loose radical extractions. 
\sa

It appears that we are now faced with the additional problems of finding out which equations
have roots with tight and/or loose  and/or natural $\ldots$-etc formulas. However, we shall soon see 
that things are not that complicated. To begin with we show that the roots of a cyclic equation 
of prime degree $p$, after the adjunction of a primitive $p^{th}$-root of unity, 
may be given a tight natural formula. In fact, to do so we need only add a Galois twist to Lagrange's 
identity 5.32.
\sas

\noindent{\bol Proposition 7.2}

{\ita Let $p$ be a prime and let $B(t)=\prod_{i=0}^{p-1}(t-\bbb_i)\in \CF[t]$ be cyclic
with Galois group $G=G_B(\CF)=\{id,\ggg,\ggg^2,\ldots ,\ggg^{p-1}\}$ where
$$
\ggg \bbb_i=\bbb_{i+1}\ess . 
\eqno 7.17
$$
Let $u$ be a primitive $p^{th}$-root of unity.  Then we have
$$
\bbb_s={1\over p}\sum_{i=0}^{p-1}\ssp u^{-is}\ssp c_i \ssp \theta^i
\bigsp (\ess for \ess\ssp s=0,..,p-1)
\eqno 7.18
$$
where each $c_i\in \CF[u]$ and $\theta$ is a root of a binomial equation
$$
t^p-\Xi=0 \bigsp (\ess \Xi\in \CF[u]\ess )
\eqno 7.19
$$
which is irreducible in $\CF[t]$.
Moreover, if $u$ is properly chosen, we have also the expansion}
$$
\theta=\bbb_o+u\bbb_1+u^2\bbb_2+\cdots +u^{p-1}\bbb_{p-1}\ess .
\eqno 7.20
$$ 
\noindent{\bol Proof}

Set 
$$
\theta_s=\sum_{i=0}^{p-1}\ssp u^{is}\ssp \bbb_i
\bigsp (\ess for \ess\ssp s=0,..,p-1)
\eqno 7.21
$$
Note that since the matrices $\|u^{rs}\|_{r,s=0,..p-1}$ and ${1\over p}\|u^{-rs}\|_{r,s=0,..p-1}$
are inverses of each other, the relation in 7.21 may be inverted to
$$
\bbb_s\ses {1\over p}\sum_{i=0}^{p-1}\ssp u^{-is}\ssp \theta_i
\bigsp (\ess for \ess\ssp s=0,..,p-1)
\eqno 7.22
$$
Note next that we may write
$$
\eqalign{
\theta_s^p&\ses \sum_{i_1,i_2,\ldots,i_p}\ssp u^{i_1+i_2+\cdots +i_p}\ssp 
\bbb_{i_1}\bbb_{i_2}\cdots \bbb_{i_p}\cr
&\ses \sum_{r=0}^{p-1}\ssp u^r\ssp \sum_{i_1+i_2+\cdots +i_p  \cong r\ess mod\ssp p}
\bbb_{i_1}\bbb_{i_2}\cdots \bbb_{i_p}\cr }\ess . 
\eqno 7.23
$$
Now it is easy to see that the polynomial
$$
P_r(\bbb_0,\bbb_1,\ldots ,\bbb_{p-1})\ses 
\sum_{i_1+i_2+\cdots +i_p \cong r\ess mod\ssp p}
\bbb_{i_1}\bbb_{i_2}\cdots \bbb_{i_p}
$$
is invariant under the cyclic shift $\bbb_i\rightarrow \bbb_{i+1}$ thus its value must necessarily
lie in $\CF$. Formula 7.23 then yields that 
$$
\theta_s^p=\Xi_s\in \CF[u]
\eqno 7.24
$$
Similarly, we can write
$$
\theta_1^{p-s}\theta_s\ses
\sum_{r=0}^{p-1}\ssp u^r\ssp \sum_{i_1+i_2+\cdots +i_{p-s}+sj \cong r\ess mod\ssp p}
\bbb_{i_1}\bbb_{i_2}\cdots \bbb_{i_{p-s}}\bbb_j
$$
and deduce from it that
$$
\theta_1^{p-s}\theta_s= d_s\in \CF[u]\ess .
\eqno 7.25
$$
Now observe that we cannot have $\theta_1=\theta_2=\cdots =\theta_{p-1}=0$
for otherwise formula 7.22 would yield $\bbb_1=\bbb_2=\cdots =\bbb_{p-1}$
which contradicts the irreducibility of $B(t)$. But this means that we can assume
$\theta_1\neq 0$. In fact, the case when  $\theta_s$ is the first one that is $\neq 0$ 
can be reduced to the case $\theta_1\neq 0$ by replacing $u$ with  $v=u^s$.
This given, we can invert $\Xi_1$ in $\CF[u]$ and, combining 7.24 with 7.25, derive that
$$
\theta_s=c_s\ssp \theta_1^s 
\eqno 7.26
$$
with
$$
c_s\ses d_s/\Xi_1\in \CF[u]\ess .
$$   
Using 7.26 in 7.22 gives 7.18 with $\theta=\theta_1$. 
Since $\theta_1$ satisfies the equation 7.19 with $\Xi=\Xi_1$ and
7.20 is simply the definition 7.21 of $\theta_1$, to complete the proof we need only show 
that the binomial
$ 
t^p-\Xi_1
$ 
is irreducible in $\CF[u]$. Now suppose it is reducible. Then Proposition 7.1 gives that it factors
in the form
$$
t^p-\Xi_1\ses (t-b)(t-ub)\cdots (t-u^{p-1})
$$
with $b\in \CF[u]$. Thus we must have $\theta_1=u^hb\in\CF[u]$ and in particular 7.18 gives
that $\bbb_1=\theta(u)$ with $\theta(t)\in\CF[t]$. Now let $H$ be the Galois group of the cyclotomic
polynomial $\PP_p=1+t+\cdots +t^{p-1}$ in $\CF[t]$ and let $H_{\bbb_1}$ denote the
Galois stabilizer of $\theta(u)$ in $H$. From the left coset decomposition
$$
H=\tau_1H_{\bbb_1}+\tau_2H_{\bbb_1}+\cdots +\tau_kH_{\bbb_1}\bigsp (\ess \tau_1=id\ess )
$$
construct the polynomial
$$
B_1(t)\ses \prod_{i=1}^k (t-\tau_i\bbb_1)\in \CF[t]\ess .
$$
Note that $B_1(t)$ must have $B(t)$ as a factor since it shares the root $\bbb_1$ with  $B(t)$
and $B(t)$ is irreducible in $\CF[t]$. Now all this leads to an impossibility since $B_1(t)$ is of
degree $k\leq p-1$. In fact, $k$ must be a divisor of the order of $H$ and $H$ is (by Corollary 7.1)
a subgroup of the cyclic group of order $p-1$. 
\sas

Thirty years before Galois and ninety years before the concept of field was introduced and 
developped by Dedekind [], Gauss, essentially showed [] that all roots of unity can be 
obtained by a sequence of {\ita normal field extensions}, in fact by {\ita tight root extractions} in our
terminology. His proof his algorithmic and is therefore
very much in style with present preferences for explicit constructions. This given, we would be amiss
not presenting it here. Gauss' argument relies on a family of remarkably beautiful identities. 
They may be stated as follows
\sa

\vbox{
\noindent{\bol Lemma 7.1} (Gauss)

{\ita Let $p$ be a prime, let $w=e^{2\pi i/p}$ and $e$ be a primitive exponent modulo $p$. For
a given factorization 
\vskip -.125truein
$$
p-1 \ses a\ssp b
\eqno 7.27
$$
set
\vskip -.125truein
$$
\Pi_a(t) \ses t+ t^{e^a}+t^{e^{2a}}+\cdots +t^{e^{(b-1)a}}\ses \sum_{i=0}^{b-1}\ssp t^{e^{ia}}\ess .
\eqno 7.28
$$
Then for any pair $0\leq r,s \leq a-1$ we have}
\vskip -.125truein
$$
\Pi_a(w^{e^r})\Pi_a(w^{e^s})\ses \sum_{i=0}^{b-1}\ssp \Pi_a(w^{e^r+e^{s+ia}})\ess .
\eqno 7.29
$$
\vskip -.125truein
\noindent {\bol Proof}
}
Note first that, for the purpose of establishing 7.28,
the polynomial $\Pi_a(t)$ is essentially invariant under the substitution $t\rightarrow t^{e^a}$.
More precisely, since $e^{p-1}\cong 1$ modulo $p$ we have
$$
\Pi_a(t^{e^a})\ses t^{e^a}+t^{e^{2a}}+\cdots +t^{e^{ba}}\ses \Pi_a(t) \bigsp (\ssp mod\ess t^p-1\ssp)
\eqno 7.30
$$
This given we have
$$
\eqalign{
\Pi_a(w^{e^r})\Pi_a(w^{e^s})&=\sum_{j=0}^{b-1}\ssp \bigl(w^{e^r}\bigr)^{e^{ja}}\ssp \Pi_a(w^{e^s})
=\sum_{j=0}^{b-1}\ssp \bigl(w^{e^r}\bigr)^{e^{ja}}\ssp \Pi_a\bigl((w^{e^s})^{e^{ja}}\bigr)\cr  
&=\sum_{j=0}^{b-1}\sum_{i=0}^{b-1} 
\ssp \bigl(w^{e^r}\bigr)^{e^{ja}}\ssp \bigl((w^{e^s})^{e^{ja}}\bigr)^{e^{ia}} 
=\sum_{i=0}^{b-1}\sum_{j=0}^{b-1} 
\ssp \bigl(w^{e^r}\bigr)^{e^{ja}}\ssp \bigl((w^{e^{s+ia}}\bigr)^{e^{ja}}\cr 
&=\sum_{i=0}^{b-1}\ssp \Pi_a(w^{e^r+e^{s+ia}})\ess .\cr}
$$
\hfill Q.E.D.
\sa

Gauss' construction of the roots of unity is best understood if illustrated
in special cases. Nevertheless it will be helpful if we start with some general remarks. 
Suppose that we want all the primitive $p^{th}$ roots of unity for a certain prime $p$. Since
they are the roots of the  cyclotomic polynomial
$$
\PP_p\ses 1 +t+t^2+\cdots +t^{p-1}\ess , 
\eqno 7.31
$$
our starting point should be the cyclicity of this polynomial with respect to the rationals. 
Following the proof of Corollary 7.2, we let $e$ be a primitive exponent modulo $p$,
let $\aaa$ be one of the roots of 7.31 and label all the roots by setting
$$
\aaa_i\ses \aaa^{e^i}\bigsp\ess\ess\ess\ess (\ess for\ess\ess i=0,1,\ldots ,p-2\ess )
\eqno 7.32
$$ 
We have seen that with this labeling the Galois group $G=G_{\PP_p}(\CQ )$ is
the cyclic group of order $p-1$ generated by the cyclic permutation
$$
\ggg \ssp \aaa_i\ses \aaa_{i+1}\ess .
\eqno 7.33
$$
Now for a given factorization $p-1=ab$ we let $G_a$ denote the cyclic subgroup of $G$ 
generated by the cycle $\ggg^a$. That is
$$
G_a=1+\ggg^a+\ggg^{2a}+\cdots +\ggg^{(b-1)a}\ess .
\eqno 7.34
$$
We then have the left coset decomposition
$$
G= G_a +\ggg G_a +\ggg^2 G_a+\cdots +\ggg^{a-1}G_a\ess .
\eqno 7.35
$$
Now we know from Theorem 6.1 that we can find a polynomial $\PS_a(\xon)\in \CQ[\xon]$
such that $\TG_\PS=G_a$. Actually in this case we have a very simple choice for $\TPS_a$, namely the
linear expression
$$
\TPS_a=\bbb_o=G_a \aaa_o=\aaa_o+\aaa_a+\aaa_{2a}+\cdots +\aaa_{(b-1)a}
\eqno 7.36
$$
Since the stability of $\TPS_a$ under $G_a$ is clear, to assure that $\TG_{\PS_a}=G_a$ we need only verify that
no other element of $G$ leaves $\TPS_a$ unchanged. But, because of 7.35, this will be so if and only if 
the conjugates
$$
\bbb_i=\ggg^i\ssp \bbb_o
\ses \aaa_i+\aaa_{i+a}+\aaa_{i+2a}+\cdots +\aaa_{i+(b-1)a}
\bigsp  (\ess for\ess\ess i=1,2,\ldots ,a-1\ess )
\eqno 7.37
$$
are all distinct. However, since
$$
\bbb_o=\sum_{r=0}^{b-1}\aaa^{e^{ra}}
\eqno 7.38
$$
the equality
$$
\bbb_i\ses \bbb_j
$$
holds true if and only if $\aaa$ is a root of the polynomial
$$
P_{ij}(t)\ses \sum_{r=0}^{b-1}t^{e^{i+ra}}\sms \sum_{r=0}^{b-1}t^{e^{j+ra}}\ess .
$$
Adding $\pm\PP_p(t)$ to $P_{ij}(t)$, if necessary to cancel the the term in  $t^{p-1}$,
we would then obtain a polynomial in $\CQ[t]$ of degree less than $p-1$ which shares a root
with $\PP_p$ contradicting the irreducibility of $\PP_p$. Thus $\TG_{\bbb_o}=G_a$
and we can use Theorem 6.5 to conclude that
\sas
\item {\ita (i)}

The polynomial
$$
B_a(t)\ses \prod_{i=0}^{a-1}(t-\bbb_i)\ess \in \ess \CQ[t]
$$
is irreducible in $\CQ[t]$ 
\sas

\item {\ita (ii)}

By adjoining $\bbb_o$ to $\CQ$ the Galois group of $\PP_p$ reduces to
$$
G_a=1+\ggg^a+\ggg^{2a}+\cdots +\ggg^{(b-1)a}\ess .
$$
\sas

\item {\ita (iii)}

The Galois group of $B_a(t)$ is the cyclic group
$$
G/G_a\ses 1+\ggg+\ggg^2+\cdots +\ggg^{a-1}\ess .
\eqno 7.39
$$
Here $\ggg$ can keep the same meaning as before since we may let it act on the $\bbb_i$'s 
as they are given by the defining identities 7.37.
\sa

This establishes that $B_a(t)$ itself is cyclic and therefore, if $a$ is a prime, we can 
construct its roots according to formula  7.18 of Proposition 7.2. All this may be
very nice but yet not very explicit! The beauty of Gauss identities is that they enable us
to compute all that we need with a minimum of effort. Indeed, 
(using the notation of Lemma 7.1), we see that we have
$$
\bbb_i\ses \Pi_a(\aaa^{e^i})\ess ,
\eqno 7.39
$$
thus we may use 7.29 to construct a multiplication table for the roots of $B_a(t)$
and obtain $B_a(t)$ itself as well as the ingredients entering in formula 7.18  quite explicitely.
\sa

This is but the first step in the construction of the roots of $\PP_p$. It reduces us to
work with $\PP_p$ in the extended field $\CQ[\bbb_o]$ which now (by the cyclicity of $B_a(t)$)
contains $\bbb_1,\bbb_2,\ldots ,\bbb_{a-1}$ as well. The next step is to factorize in turn 
the new Galois group of $\PP_p$, which as we have seen reduces to
$$
G_a=1+\ggg^a+\ggg^{2a}+\cdots +\ggg^{(b-1)a}\ess ,
$$
and then proceed to split each of the $\bbb_i$ into sums $\la_j$ of powers of $\aaa$ which are invariant
under a normal subgroup of $G_a$. Then, after the adjonction of the $\la_j$ to $\CQ[\bbb_o]$,
we further reduce the Galois group of $\PP_p$ to this normal subgroup. We proceed in
this manner until the $\bbb_i$'s are split all the way down to their individual summands,
which are of course the roots of $\PP_p$. To describe the process explicitely at this point and
in full generality would require more notation that would only blur the beauty of the argument.
Imitating established tradition, we will avoid this difficulty  by just saying that after we split the number $p-1$
into its prime factors. 
$$
p-1=p_1p_2 \cdots p_{k_o}\ess\ess\ess {\rm with}\ess\ess\ess p_1\geq p_2\geq \cdots \geq p_{k_o}\ess ,
\eqno 7.40
$$
we proceed to construct the {\ita composition series} of $G=1+\ggg+\cdots +\ggg^{p-1}$ 
\def \norr {\triangleright}
$$
G\norr_{p_1}G_1\norr_{p_2}G_2\norr_{p_3}\cdots \norr_{p_{k_o}}G_{k_o}=\{id\}\ess .
$$
Then by a sequence of adjunctions $\xi_k\in \CQ[\aaa_o,\aaa_1,\ldots, \aaa_{p-1}]$, with
$\TG_{\xi_k}=G_k$ we arrive at a final tight and natural expression for the roots of $\PP_p$.
We must also point out that this argument can only be completed by an induction process.
Since at each step, as we have seen we need to adjoin primitive roots of unity of lower order
which inductively must be assumed to have already been given tight formulas. Before we indulge
into this type of mental gymnastics it will be good to work out a few revealing examples.
\sa

\item{$\underline {p=13}$:}

Here, $p-1=12$ and we may take $e=2$ as a primitive exponent $mod\ess 13$. 
This given,  from Corollary 7.2 we get that the Galois group of the equation
$$
\PP_{13}(t)=1+t+t^2+\cdots + t^{12}
$$
with respect to the field $\CQ$ of rational numbers,
is the cyclic group generated by the operation  $t\rightarrow t^2$. 
More precisely, if $\aaa$ is our desired $13^{th}$ root of unity and the roots of $\PP_{13}(t)$
are labeled by setting
$$
\aaa_o=\aaa\scs \aaa_1=\aaa^2\scs \aaa_2=\aaa^{2^2}\scs \ldots , \aaa_{11}=\aaa^{2^{11}}\ess ,
$$
then $G=G_{\PP_{13}}(\CQ)$ is generated by the cyclic permutation
$$
\ggg \aaa_i \ses \aaa_i^2 \ses \aaa_{i+1}\ess .
$$
This gives
$$
\matrix{
\aaa_o&\ses &\aaa^{2^0}&=&w \cr
\aaa_1&\ses&\aaa^{2^1}&=&\aaa^{2}\cr
\aaa_2&\ses&\aaa^{2^2}&=&\aaa^{4}\cr
\aaa_3&\ses&\aaa^{2^3}&=&\aaa^{8}\cr
\aaa_4&\ses&\aaa^{2^4}&=&\aaa^{3}\cr
\aaa_5&\ses&\aaa^{2^5}&=&\aaa^{6}\cr
\aaa_6&\ses&\aaa^{2^6}&=&\aaa^{12}\cr
\aaa_7&\ses&\aaa^{2^7}&=&\aaa^{11}\cr
\aaa_8&\ses&\aaa^{2^8}&=&\aaa^{9}\cr
\aaa_9&\ses&\aaa^{2^9}&=&\aaa^{5}\cr
\aaa_{10}&\ses&\aaa^{2^10}&=&\aaa^{10}\cr
\aaa_{11}&\ses&\aaa^{2^11}&=&\aaa^{7}\cr}
\eqno 7.41
$$
We start by factoring the Galois group of $\PP_{13}$ 
$$
1+\ggg+\cdots+\ggg^{11}\ses 
1+\ggg^3+\ggg^6+\ggg^{9}+
\ggg (1+\ggg^3+\ggg^6+\ggg^{9})+
\ggg^2(1+\ggg^3+\ggg^6+\ggg^{9})\ess .
\eqno 7.42
$$
Then we seek for an element $\bbb_o\in\CQ[\aaa_o\aaa_1,\ldots,\aaa_{11}]$ whose Galois stabilizer is 
$$
\TG_{\bbb_o}\ses 1+\ggg^3+\ggg^6+\ggg^{9}
$$
We may take
$$
\bbb_o\ses \aaa_o+\aaa_3+\aaa_6+\aaa_9=\aaa+\aaa^8+\aaa^{12}+\aaa^5\ses \Pi_3(\aaa)\ess .
\eqno 7.43
$$
Its conjugates are 
$$
\bbb_1\ses \aaa_1+\aaa_4+\aaa_7+\aaa_{10}=\aaa^2+\aaa^3+\aaa^{11}+\aaa^{10}\ses\Pi_3(\aaa_1)\ess .
\eqno 7.44
$$
and
$$
\bbb_2\ses \aaa_2+\aaa_5+\aaa_8+\aaa_{11}=\aaa^4+\aaa^6+\aaa^9+\aaa^7\ses\Pi_3(\aaa_2)\ess .
\eqno 7.45
$$
Now, using the Gauss relations with parameters $p=13$, $a=3$, $b=4$ and $e=2$ we get the 
multiplication table
$$
\matrix{
\bbb_o\bbb_o=4+\bbb_1+2\bbb_2 & \bbb_o\bbb_1=\bbb_o+2\bbb_1+\bbb_2 & \bbb_o\bbb_2=2\bbb_o+\bbb_1+\bbb_2 \cr
\bbb_1\bbb_o=\bbb_o+2\bbb_1+\bbb_2 &\bbb_1\bbb_1=4+2\bbb_o+\bbb_2&\bbb_1\bbb_2=\bbb_o+\bbb_1+2\bbb_2 \cr
\bbb_2\bbb_0=2\bbb_o+\bbb_1+\bbb_2 &\bbb_2\bbb_1=\bbb_o+\bbb_1+2\bbb_2 & \bbb_2\bbb_2=4+\bbb_o+2\bbb_1 \cr
}
\eqno 7.46
$$
From which we derive that 
$$
(t-\bbb_o)(t-\bbb_1)(t-\bbb_2)\ses 1-4t+t^2+t^3\ess .
$$
Theorem 6.5 then gives that this equation is cyclic with Galois group 
$$
1+\ggg+\ggg^2
$$
So if we let $u$ denote a primitive cube root of unity, that is
$$
u={-1+{\root 2\of {-3}}\over 2}
$$
and set
$$
\eqalign{
\theta_o&\ses \bbb_o+\bbb_1+\bbb_2\ess ,\cr
\theta_1&\ses \bbb_o+u\bbb_1+u^2\bbb_2\ess ,\cr
\theta_2&\ses \bbb_o+u^2\bbb_1+u^4\bbb_2 \ess , \cr}
$$
then formula 7.22 gives
$$
\bbb_o=(\theta_o+\theta_1+\theta_2)/3\ess .
$$
Note that we must have $-\theta_o=1$ since it must equal the coefficient of $t^{11}$ in $\PP_{13}$.
To get $\bbb_o$ into the form given in 7.18 we must compute the coefficient
$$
c_2\ses \theta_2/\theta_1^2\ess .
$$
Now we can write 
$$
c_2\ses {\theta_2\theta_1\over \theta_1^3}\ess .
$$
and the table in 7.46 gives
$$
\theta_2\theta_1\ses 13
\ess\ess\ess ,\ess\ess\ess
\theta_1^3\ses -13(4+3u)\ses {13^2\over -1+3u}
$$
Thus 
$$
c_2\ses (-1+3u)/13
$$ 
and 
$$
\bbb_o\ses {-1+\theta_1+(-1+3u)\theta_1^2/13 \over 3}\ess .
\eqno 7.47
$$
where we may write
$$
\theta_1\ses {\root 3 \of {-13(4+3u)}}\ess ,
$$
which is assured by Proposition 7.2 to be a tight radical in $\CQ[u]$. 
\sas

\noindent
By Theorem 6.5 the Galois group of $\PP_{13}$ in $\CQ[\bbb_o]$ is
$$
1+\ggg^3+\ggg^6+\ggg^9
$$
So our next step is to construct a polynomial in $\aaa_o,\aaa_1,\ldots ,\aaa_{11}$ whose
stabilizer in $1+\ggg^6$. We may take
$$
\la_o\ses \aaa_o +\aaa_6 \ses \aaa +\aaa^{2^6}\ses \aaa+\aaa^{-1}
$$
Now we need to use the Gauss relations with $a=6$ $e=2$. In fact, if we set
$$
\la_r\ses \Pi_6(\aaa^{2^r})\bigsp (\ess r=0,1,\ldots ,5\ess )
$$
then 
$$
\la_o\ses \aaa_o+\aaa_6
\ess\ess\ess{\rm and} \ess\ess\ess
\la_3\ses \aaa_3+\aaa_9
$$
thus
$$
\bbb_o\ses \la_o+\la_3\ess .
$$
From Gauss formula (or even by direct computation in this case) we get that
$$
\la_o\times \la_3\ses \la_2+\la_5\ses \bbb_2
$$
and $\la_o$ can be obtained by solving the equation
$$
(t-\la_o)(t-\la_3)\ses t^2-\bbb_o t+\bbb_2\ess .
$$
This gives that
$$
\aaa+\aaa^{-1}=\la_o\ses {-\bbb_o+{\root 2 \of {\bbb_o^2 -4\bbb_2}}\over 2}
\eqno 7.48
$$
If we prefer to write $\la_o$ only in terms of $\bbb_o$ the we resort again to
the table in 7.46 which gives 
$$
\bbb_1=-\bbb_o^2-2\bbb_o+2
\ess\ess\ess{\rm and} \ess\ess\ess
\bbb_2=\bbb_o^2+ \bbb_o-3
$$
this agrees with the fact (implied by the cyclicity of the equation satisfied by $\bbb_o$)
that both $\bbb_1$ and $\bbb_2$ must be in $\CQ[\bbb_o]$.
We should note again that, as assured by Proposition 7.2, the square root in 7.48 must necessarily
be tight in $\CQ[\bbb_o]$. 
\sas

The last step is to reduce the Galois group of $\PP_{13}$ to the identity by adjoining
$\aaa$. From 7.48 we get that $\aaa$ is obtained by solving the equation
$$
(t-\aaa)(t-\aaa^{-1})\ses t^2-\la_ot +1=0
$$
which gives that
$$
\aaa\ses {\la_o+{\root 2 \of {\la_o^2-4}}\over 2}
\eqno 7.49
$$
In summary, by combining 7.47, 7.48 and 7.49, we can construct a tight natural 
formula for all the primitive $13^{th}$ roots of unity. In fact the successive 
adjunctions may be taken to be
$$
\bbb_o=\aaa_o+\aaa_3+\aaa_6+\aaa_9
\ess\ess ,\ess\ess
\la_o\ses \aaa_o+\aaa_6
\ess\ess ,\ess\ess
\aaa_o=\aaa
$$
The corresponding reductions of the Galois group of $\PP_{13}$ being
$$
1+\ggg+\cdots+\ggg^{11}\ess \norr_3\ess 1+\ggg^3+\ggg^6+\ggg^9\ess\norr_2\ess 1+\ggg^6\ess\norr_2\ess \{id\}
$$
\sas

\noindent{\bol Remark 7.1}

The fact that in the previous calculation we found that $\theta_1\theta_2=13\in \CQ$ is not an 
accident. In general, for any cyclic equation, formulas 7.21 give that
$$
f\ses \theta_1\theta_2\cdots \theta_{p-1}\ses
\sum_{r=0}^{p-1}u^r 
\ssp\sum_{i_1+2i_2+\cdots +(p-1)i_{p-1}\cong r \ssp mod \ssp p}
\bbb_{i_1}\bbb_{i_2}\cdots \bbb_{i_{p-1}}\ess .
$$
Since the coefficient of $u^r$ is clearly invariant under the cyclic shift $\bbb_i\rightarrow \bbb_{i+1}$,
we should not be surprised if $f$ comes out to be in $\CQ[u]$. However, it is easy to see that $f$ also
remains unchanged by the replacement of $u$ in $f$ by any other primitive $p^{th}$-root of unity.
This forces $f\in \CQ$ as well.
\sas

\item {$\underline {p=11}$:}

Here we may take $e=2$ again and set $\aaa_i=\aaa^{2^i}$ for $i=0,1,\ldots ,9$ with $\aaa$
our desired primitive $11^{th}$-root of unity. Letting  $\ggg\aaa_i=\aaa_{i+1}$ again 
Corollary 7.2 gives 
$$
G_{\PP_{11}}(\CQ)\ses
1+\ggg+\ggg^2+\cdots +\ggg^9
$$
Since $p-1=10=5\times 2$ we can at once reduce this Galois group to the subgroup
$$
1+\ggg^5
$$
by the adjonction of
$$
\bbb_o\ses \aaa_o+\aaa_5\ess .
\eqno 7.50
$$
To use Gauss machinery with $p=11$ $e=2$ and $a=5$ we set
$$
\bbb_r\ses \Pi_5(\aaa_r)\ess . \bigsp\bigsp (\ess r=0,1,\ldots ,4\ess )
$$
In this case a repetitive application of 7.29 gives
$$
(t-\bbb_o)(t-\bbb_1)\cdots (t-\bbb_4)\ses 1+3t-3t^2-4t^3+t^4+t^5\ess ,
$$
which by Theorem 6.4 must be cyclic with Galois group
$$
1+\ggg+\ggg^2+\ggg^3 +\ggg^4
$$
We can then solve it with the formulas of Proposition 7.2. So we pick
a primitive $5^{th}$-root of unity $u$ and set
$$
\theta_r\ses \sum_{i=0}^4 \ssp u^{ri}\ssp \bbb_i
\bigsp (\ess for\ess r=0,1,\ldots ,4\ess )\ess .
$$
Using the Gauss identities, (this time MAPLE comes in handy) we get that
$$
\theta_1^5\ses -11u(26+20u-15u^2 +10u^3)\ess .
$$
So $\theta_1$ is obtained by the  tight root extraction
$$
\theta_1\ses {\root 5 \of{ -11u(26+20u-15u^2 +10u^3)}}
\eqno 7.51
$$
and $\bbb_o$ is then given by
$$
\bbb_o\ses (\theta_o+\theta_1+\cdots +\theta_4)/5
$$
Now again we have $\theta_o=-1$. To express $\bbb_o$ entirely in terms of $\theta_1$
we need to compute the ratios $c_i=\theta_i/\theta_1^i$. We shall compute $c_2$ here
and leave it to the reader to compute the other $c_i$'s. 

In complete agreement with our Remark 7.1 we find that
$$
\theta_1\theta_2\theta_3\theta_4\ses 121=11^2\in\CQ
$$
So we may write
$$
\theta_2=c_2\theta_1^2\ses {121\over \theta_1^3\theta_3\theta_4}\ssp \theta_1^2
$$
Again the Gauss formulas (and MAPLE) give that
$$
\theta_1^3\theta_3\theta_4\ses 121u(2-2u-u^2)\ess ,
$$
and we can then easily get that
$$
\theta_2\ses {4+2u+2u^2+u^3\over 11 u}\ess \theta_1^2\ess .
$$
Since $2^5\cong -1$ mod $11$ we see that
$$
\bbb_o\ses \aaa+\aaa^{-1}
$$
so in one more step we can get our desired $\aaa$ by solving the equation
$$
(t-\aaa_o)(t-\aaa_5)\ses t^2-\bbb_ot+1\ses 0\ess .
$$
This gives again
$$
\aaa\ses {\bbb_o+{\root 2 \of{\bbb_o-4}}\over 2}
$$
In accordance with the fact that $p-1=5\times 2$ here we have only needed two natural adjunctions
to reduce the Galois group to the identity. Namely, $\aaa_o+\aaa_5$ and $\aaa_o$.  
\sas

\item {$\underline {p=7}$:}

We will be brief here. In this case we must take $e=3$ as a primitive exponent and
set $\aaa_i=\aaa^{3^i}$ for $i=0,1,\ldots ,5$, where $\aaa$ is our desired $7^{th}$-root of unity.
The Galois group in this case is
$$
G_{\PP_7}(\CQ)\ses 1+\ggg+\cdots +\ggg^5\bigsp (\ess with \ess \ggg \aaa_i=\aaa_{i+1}\ess )
$$
Since
$p-1=3\times 2$ we start with
$$
\bbb_r\ses \Pi_3(\aaa_r)\bigsp (\ess for \ess r=0,1,2\ess )
$$
Then Gauss' formulas yield us the table
$$
\matrix{
\bbb_o\bbb_o=2+\bbb_2 & \bbb_o\bbb_1=\bbb_1+\bbb_2 &\bbb_o\bbb_2=\bbb_o+\bbb_1 \cr
\bbb_1\bbb_o=\bbb_1+\bbb_2 & \bbb_1\bbb_1=2+\bbb_o &\bbb_1\bbb_2=\bbb_o+\bbb_2 \cr
\bbb_2\bbb_o=\bbb_o+\bbb_1 & \bbb_2\bbb_1=\bbb_o+\bbb_2 &\bbb_2\bbb_2=2+\bbb_1 \cr}
$$
From which we get that
$$
(t-\bbb_o)(t-\bbb_1)(t-\bbb_2)\ses t^3+t^2-2t-1\ess .
$$
Theorem 6.5 gives that this polynomial is cyclic with Galois group $1+\ggg+\ggg^2$. 
So we may use Proposition 7.2 with $p=3$ and set again
$$
\eqalign{
\theta_o&\ses \bbb_o+\bbb_1+\bbb_2\ess ,\cr
\theta_1&\ses \bbb_o+u\bbb_1+u^2\bbb_2\ess ,\cr
\theta_2&\ses \bbb_o+u^2\bbb_1+u^4\bbb_2 \ess , \cr}
$$
This gives us
$$
\bbb_o=(-1+\theta_1+\theta_2)/2
$$
Using the table  we easily derive that
$$
\theta_1\theta_2\ses 7\ess\ess\ess and \ess\ess\ess \theta_1^3\ses -7-21 u
$$
so we may write
$$
\theta_2\ses {7\over \theta_1}\ses  {7\over \theta_1^3}\ess \theta_1^2\ses 
{2+3u\over 7}\ess \theta_1^2\ess .
$$
So we get
$$
\theta_1\ses {\root 3 \of {-7-21u}}
$$
and
$$
\bbb_o={-1+\theta_1+ (2+3u)\theta_1^2/7\over 2}\ess .
$$
Now here 
$$
\bbb_o=\aaa_o+\aaa_3\scs 
\bbb_1=\aaa_1+\aaa_4\scs 
\bbb_2=\aaa_2+\aaa_5
$$
and since $\aaa_3=\aaa^{-1}$ we can find $\aaa$ by solving the equation
$$
(t-\aaa_o)(t-\aaa_3)\ses t^2-\bbb_ot+1
$$
which gives 
$$
\aaa={\bbb_o+{\root 2\of {\bbb_o^2-4}}\over 2}
$$
 
\sas

\item {$\underline {p=5}$:}

Here we may take the exponent $e=2$. So if $\aaa$ denotes our desired primitive
$5^{th}$-root, we need to set $\aaa_i=\aaa^{2^i}$ for $i=0,1,\ldots ,3$.
Since $p-1=2\times 2$ we start by constructing the two
elements
$$
\bbb_o=\aaa_o+\aaa_2\ess\ess\ess ,\ess\ess\ess \bbb_1=\aaa_1+\aaa_3
$$
by solving the quadratic
$$
(t-\bbb_o)(t-\bbb_1)\ses t^2+t-1\ess .
$$
Thus we may take
$$
\bbb_o\ses {-1+{\root 2 \of 5}\over 2}\ess\ess\ess,\ess\ess\ess
\bbb_1\ses {-1-{\root 2 \of 5}\over 2}
\ess ,
$$
and $\aaa$ can be obtained by solving
$$
(t-\aaa_o)(t-\aaa_2)\ses t^2-\bbb_ot+1\ses 0\ess .
$$
Since $\bbb_o^2=\bbb_1+2$ we finally obtain
$$
\aaa={\bbb_o+{\root 2 \of {\bbb_1-4}}\over 2}
\ses 
{
 {
   {-1+{\root 2 \of 5}\over 2}
   + i\ssp
   {\root 2\of 
              {
                5+{\root 2 \of 5}\over 2
              }
   }
 }
\over 2
}
$$
\sa

We are now in a position to establish the basic result of this section
\sas

\noindent{\bol Theorem 7.2}

{\ita Every cyclic equation can be solved by a sequence of tight and natural
radical extractions}
\sas

\noindent{\bol Proof}

Suppose that $\TE_n(t)$ is cyclic in $\CF$, and let $\aon$ be the labeling of
its roots under which
$$
G_{\TE_n}(\CF)\ses 1+\ggg+\ggg^2+\cdots +\ggg^{n-1}\ess .
\bigsp \bigsp(\ssp with \ess  \ggg\aaa_i=\aaa_{i+1} \ess )
$$
If $p_1\geq p_2\geq \cdots \geq p_m$ are the primes factoring $n$ we proceed by
constructing a sequence of subgroups $G_k$ yielding the composition series
$$
G_{\TE_n}(\CF)\norr_{p_1}G_1\norr_{p_2}G_2\norr_{p_2}\cdots \norr_{p_m}G_m=\{id\}
$$
This given we construct (via Theorem 6.1) a sequence of polynomials
$\PS_k(\xon)\in \CF[\xon]$ with $\TG_{\PS_k}=G_k$. Now Theorem 6.4 and
$\TG_{\PS_k}\nor_{p_k} \TG_{\PS_{k-1}} $ give that 
$$
\CF[\TPS_k]\son \CF[\TPS_{k-1}]\ess .
$$
Moreover from Theorem 6.5 we also deduce that $\TPS_k$ is a root of a polynomial
$B_{\TPS_k}(t)\in \CF[\TPS_{k-1}][t]$ of degree $p_k$ which is cyclic in $\CF[\TPS_{k-1}]$.
Thus its roots may be constructed by means of the formulas given by Proposition 7.2.
Let $\root p_k \of \Xi_k$ with $\Xi_k\in \CF_{k-1}$ denote the 
radical we must extract to obtain $\TPS_k$. Now Proposition 7.2 assures that  
$\root p_k \of \Xi_k$ is tight, and formula 7.20 gives that this a natural radical
extraction as long as we are in possession of a primitive $p_k^{th}$ root of unity.
Finally, the identity in 7.18 shows that  $\CF[\TPS_k]$ may also be obtained from $\CF[\TPS_{k-1}]$ 
by the adjonction of $\root p_k \of \Xi_k$ itself. This given, since the Galois group
of $\TE_n$ reduces to the identity after the extraction of the last radical $\root p_m \of \Xi_m$,
the process will yield a tight natural formula for each of the roots of $\TE_n$.

We should note that the tight extraction of prime $p^{th}$  roots of unity can also
included in this process. This is because, as we have seen, the formulas giving a primitive
$p_k^{th}$-root depend on the solution of cyclic equations whose degrees are prime factors
of $p_k-1$. But since each $p_k\leq n$, we see that we are appropriately setup for an induction
argument. We can in fact assume from the start that the Theorem is true for all cyclic equations 
of prime degree less or equal than a certain prime $p$ then carry out the constrution
outlined above for all cylic equations of any degree $n$ whose prime factors are all
less or equal to the next prime. The induction can of course start with $p=2$ where
the Theorem is easily verified to be true.
\sa
\sa

\centerline{\bol Bibliography}

\sa

\item {[1]}  Elwyn R. Berlekamp,
{\ita ALGEBRAIC CODING THEORY},
 Revised 1984 Edition,
 Aegean Park Press. 1-478
 
\sas

\item {[2]}  Edgar Dehn,
{\ita ALGEBRAIC EQUATIONS},
 An Introduction to the Theories of Lagrange and Galois,
 Dover Publications. 1-208
 
\sas

\sa

The author of the lecture notes has used the contents of these books in a different order. Inspired by the contents of these two books and guided by their contents. The author just used the results without using their proofs but guided by the results of Galois Theory as obtained by Galois himself.

\end

\def\picture #1 by #2 (#3){
  \vbox to #2{
    \hrule width #1 height 0pt depth 0pt
    \vfill
    \special{picture #3} 
    }
  }

\def\scaledpicture #1 by #2 (#3 scaled #4){{
  \dimen0=#1 \dimen1=#2
  \divide\dimen0 by 1000 \multiply\dimen0 by #4
  \divide\dimen1 by 1000 \multiply\dimen1 by #4
  \picture \dimen0 by \dimen1 (#3 scaled #4)}
  }